\RequirePackage{lineno}
\documentclass[12pt,pra,aps,amssymb,amsfonts,amsmath,tightenlines]{revtex4}%aps instead of rmp
\usepackage{placeins}
\usepackage{dsfont}
\usepackage{graphicx}
\usepackage{caption}
\usepackage{amssymb,amsfonts,amsthm}
\usepackage{color}
\usepackage{verbatim}
\usepackage[normalem]{ulem}
\usepackage{enumerate}
\usepackage{mathrsfs}
\usepackage{bm,float}
\usepackage{mathtools}
\usepackage{extarrows}

\usepackage{textgreek}
\usepackage{lgreek}
\usepackage{gensymb}

\newtheorem{Proposition}{Proposition}

\usepackage{multirow}
\usepackage{rotating}
\usepackage{caption}
\usepackage{booktabs,siunitx,array,threeparttable}
\usepackage{tikz} 
\usetikzlibrary{shapes.geometric,calc,shapes,intersections,patterns,angles,quotes,through,backgrounds,decorations.markings,graphs}
\tikzset{myDot/.style={draw,shape=circle,fill=black,minimum size=6.0pt,inner sep=0.0pt,line width=1pt}} %controls the features of the vertices
\tikzset{myDot2/.style={draw,shape=circle,fill=white,minimum size=6.0pt,inner sep=0.0pt,line width=1pt}}
\tikzset{edgeLabel/.style={font=\scriptsize, red!95!black,right}} %color and boldness of chord labels; add "sloped" to get labels that run along the edges
\tikzset{extendLine/.style 2 args={shorten >=-{#2},shorten <=-{#1}}} %requires parameters saying how far to extend edges left and right beyond vertices
\tikzset{myLine/.style={line width=.3mm}} %general lines in the figures
\tikzset{myDottedLine/.style={dotted, line width=.3mm}} %dotted lines
\tikzset{myDashedLine/.style={dashed, line width=.3mm}} %dashed lines
\tikzset{myBoldLine/.style={line width=.6mm}} %bold lines, if we need them
\begin{document}
\title{Configurations, Tessellations and Tone Networks}
\author{Jeffrey~R.~Boland$^1$ and Lane~P.~Hughston$^{2,3}$}

\affiliation{$^1$Syndikat LLC, 215 South Santa Fe Avenue, Los Angeles, California 90012, USA\\
$^{2}$School of Computing, Goldsmiths University of London, New Cross, London SE14\,6NW, UK\\
$^{3}$Artificial Intelligence and Mathematics Research Lab, James Carter Road, Mildenhall, Bury St Edmunds IP28\,7DE, UK}
%\date{\today}
%ABSTRACT HERE
\begin{abstract}
\noindent 
The tonnetz, which is commonly represented as a tessellation of the plane by a triangular network of tones, can also be represented as a bipartite graph of degree three with twelve vertices denoting major triads and twelve vertices denoting minor triads. We show that this Levi graph can be realized geometrically as a system of twelve points and twelve lines in $\mathbb R^2$ with the property that three points lie on each line and three lines pass through each point, in a configuration $\{12_3\}$ of Daublebsky von Sterneck type D222. This tonnetz configuration, alongside various generalizations thereof, can be used as a new basis for the composition and analysis of music. 
\\
%KEY WORDS HERE
\begin{center}
{\scriptsize {\bf Keywords: 
Music and mathematics, tonnetz, configuration, Levi graph, tessellation, Tristan genus, \\dominant seventh, half-diminished seventh, Tristan und Isolde,  G\"otterd\"ammerung, Parsifal} }
\end{center}
\vspace{0.05cm}
%\keywords{Key Words}
%\subclass{MSC code\and JEL classification code}
\end{abstract}

\maketitle
%Section I
%
%LINE NUMBERS
%\linenumbers
%LINE NUMBERS
%%%%%%
%SECTION I
\section{Introduction}
\label{sec:Introduction}
%%%%%%%

\begin{figure}[htbp] %triangular tessellation
\centering
\includegraphics[clip,scale=0.40]{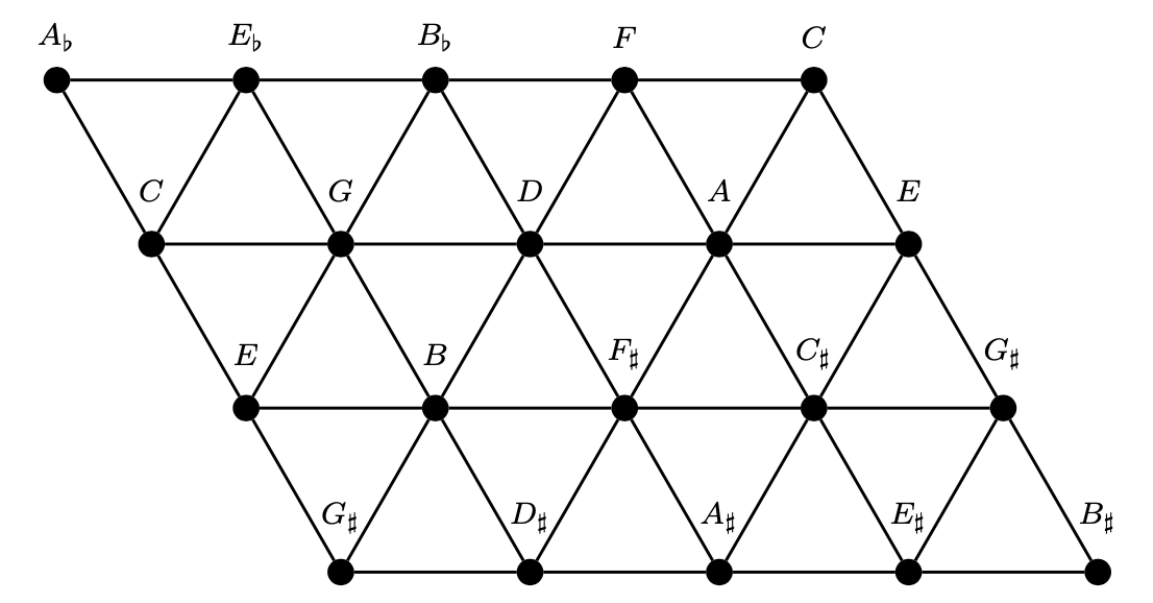}
\vspace{0.50cm}
  \caption{The tonnetz as a tessellation of the Euclidean plane by regular triangles, with notes of the chromatic scale situated at the vertices. Moving from left to right along the horizontal axis, we see that notes ascend in fifths.  Moving upwards to the right at  60$\degree$ from the horizontal axis,  notes ascend in minor thirds.  Moving downwards to the right at  60$\degree$,  notes ascend in major thirds.  }
\label{Triangular_tessellation}
\end{figure}

\begin{figure}[htbp] %chicken-wire tessellation
\centering
\includegraphics[clip,scale=0.40]{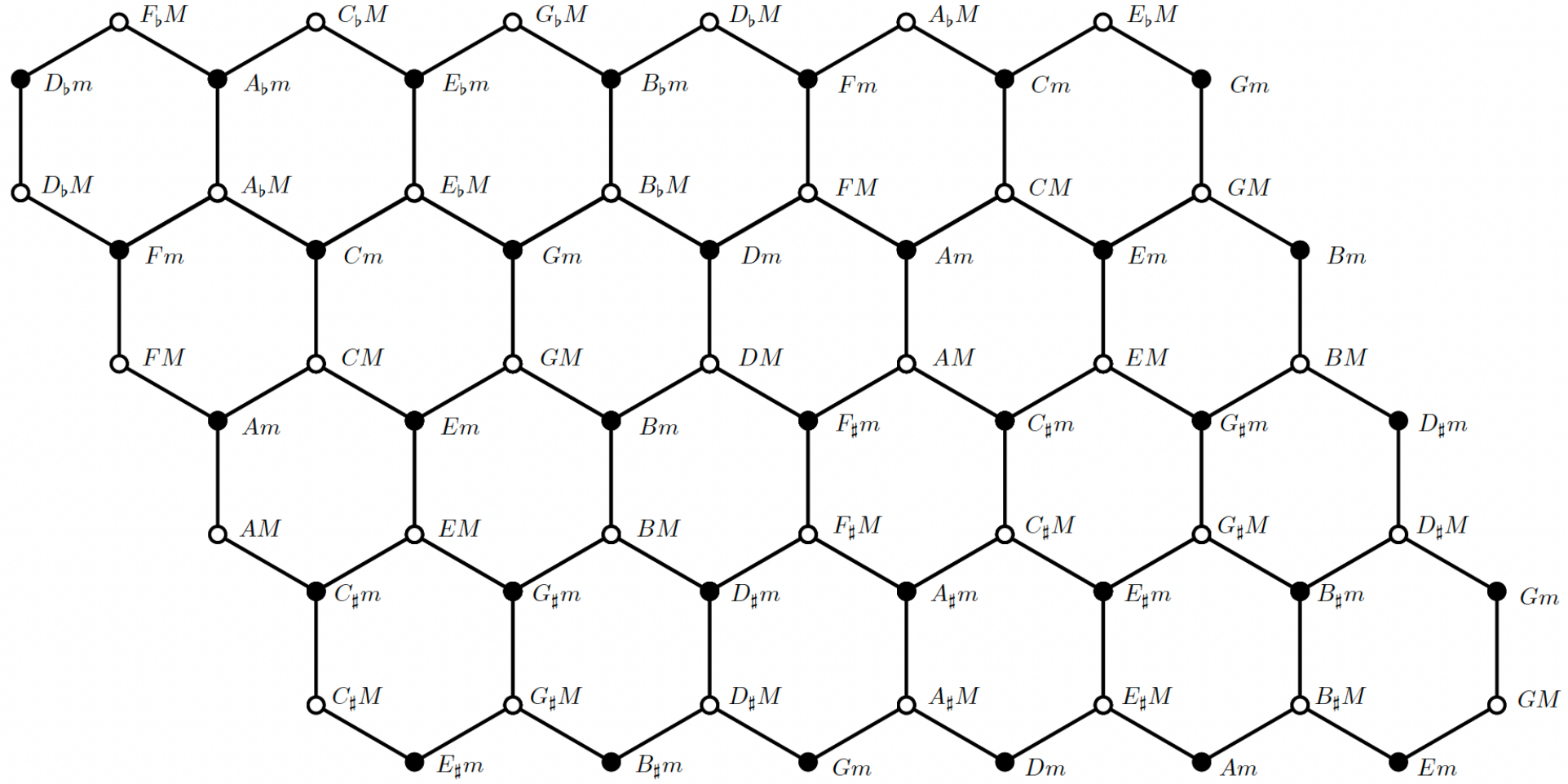}
\vspace{0.50cm}
  \caption{The tonnetz as a tessellation of the Euclidean plane by regular hexagons.  Each white vertex represents a major triad and each black vertex represents a minor triad. Three minor triads join each major triad and three major triads join each minor triad. The vertical edges such as $\{CM, Cm\}$ represent {\it parallel} ($\rm \bf P$) transformations.  The edges running from left to right upward 30$\degree$  from the horizontal such as $\{Am, CM\}$ represent {\it relative} ($\rm \bf R$) transformations. The edges running from left to right downwards 30$\degree$  from the horizontal such as $\{CM, Em\}$ represent {\it Leittonwechsel} ($\rm \bf L$) transformations.}
\label{tonnetz_tessellation}
\end{figure}

\noindent The tonnetz,  which consists of an intricate network of relations between the twelve notes of the chromatic scale, was originally introduced by the mathematician Leonhard Euler (1739) and later developed in works of Ernst Naumann (1858), Arthur von Oettingen (1866), Hugo Riemann (1880), and others \footnote{We treat the word ``tonnetz'' spelled with a lower case ``t'' and its plural ``tonnetze'' as words of English.}.  Investigation of the tonnetz from a modern standpoint in the setting of equal temperament was carried forward in works of Waller (1978), Lewin (1982, 1987), Hyer (1989), Cohn (1996, 1997, 1998), Douthett and Steinbach (1998), and others. The tonnetz is usually presented as a network of tones, in the form of a tessellation of $\mathbb R^2$ by equilateral triangles, as shown in Figure \ref{Triangular_tessellation}, with a note of the chromatic scale at each vertex.  Figure \ref{tonnetz_tessellation} shows an alternative representation, dual to that of Figure  \ref{Triangular_tessellation}, as a tessellation of the plane by regular hexagons, in which the vertices are consonant triads.  Both versions of the tonnetz, viewed as infinite graphs, are familiar to music theorists. See, for example, Tymoczko (2011) at page 413. The chicken-wire representation of Figure \ref{tonnetz_tessellation} can be found in Douthett and Steinbach (1998) at page  248.  But the tessellations  in Figures  \ref{Triangular_tessellation} and \ref{tonnetz_tessellation} are more than graphs:  they are geometric objects in the Euclidean plane. 
Here we propose to introduce another representation of the tonnetz as a geometric object. It will be shown that the tonnetz can be visualized as a configuration of twelve lines and twelve points in $\mathbb{R}^2$, with the property that three of the points lie on each line and three of the lines pass through each point. 
The construction proceeds in several steps.
First, we interpret the tonnetz as a finite bipartite graph of degree three and girth six with the two types of vertices corresponding to major triads and minor triads, respectively, and the edges encoding adjacency within the tonnetz.  
Then we realize this Levi graph as an incidence geometry consisting of twelve points (major triads) and twelve lines (minor triads) in $\mathbb R^2$, in the form of a configuration $\{12_3\}$. 

Finally, we identify this configuration as the self-dual $\{12_3\}$ of type D222 introduced by Daublebsky von Sterneck (1895).  The tonnetz is widely regarded as essential to the modern understanding of music based on triadic harmony -- so it is both unexpected and important to see such a connection between the tonnetz and a basic construction in Euclidean geometry.  
\vspace{-0.20cm}

%%%%%%%%%
%%%%%

Mathematicians have been studying configurations for centuries and there is a vast literature on the subject. A configuration $\{m_r, n_k\}$ is a collection of $m$ points and $n$ lines such that $r$ of the $n$ lines are incident with each point and $k$ of the $m$ points lie on each line.  If the number of points equals the number of lines and the number of points lying on each line equals the number of lines passing through each point, then the configuration is ``balanced'', in the terminology of Gr\"unbaum (2009). The abbreviation $\{n_k\}$ is used to denote a balanced configuration of type $\{n_k, n_k \}$.

Our central thesis is that configurations give rise to abstract systems of musical relations, in which chords of different types, such as the major and minor triads of the tonnetz, can be represented in a mathematically precise way by geometric elements such as points and lines in the Euclidean plane, where various specific chord relation categories, efficient voice-leading being amongst them, are modelled by adjacency relations among the corresponding geometric elements.  
We shall be mainly concerned with balanced configurations and we point to examples of new music systems.

%%%%%%%%%
%%%%%

We presume some familiarity with the Eulerian tonnetz, and indeed when we refer to ``the tonnetz'' without further qualification, we usually mean the Eulerian tonnetz in its modern form based on equal temperament. To establish some terminology and notation we develop the relevant ideas {\it ab initio} and we present a brief construction of the tonnetz that only involves elementary ideas.   

We fix the chromatic scale  and give the notes their usual names $C, C_{\sharp}/D_{\flat}, D, D_{\sharp}/E_{\flat}$, and so on, making no distinction between enharmonic pairs.  We use curly brackets to denote an unordered set of notes or chords. Square brackets denote an ordered set or sequence of notes. Round brackets are used for sequences of chords. Angled brackets are used in the case of a cycle of distinct chords, in which case the final written chord repeats the initially written one, to signal that the sequence closes. 
We use the term dyad for an unordered pair of distinct notes, for example $\{C,E\}$. We write 
$\{C, E, G\}$ to represent an unordered triad of the notes $C$, $E$, $G$, each of which might be located within any octave. In such a context, $C$ denotes an equivalence class of notes modulo octaves. Then we write $CM = \{C, E, G\}$ for the ``abstract'' $C$-major triad (that is, the triad of pitch classes).   

When we write $CM = [C, E, G]$ for an ordered triad in root position, it is usually with a specific choice of pitches in mind for the three notes, dictated by context. Then we write $CM^{(1)} = [E, G, C]$ and $CM^{(2)} = [G, C, E]$. Similarly, we write $Cm = [C, E_{\flat}, G]$, $Cm^{(1)} = [E_{\flat}, G, C]$ and $Cm^{(2)} = [G, C, E_{\flat}]$.
In a situation where the inversions are important, we write $CM^{(0)} = [C, E, G]$ for the triad in its root position. 
There are three further permutations of a triad, $[E, C, G]$, $[G, E, C]$ and $[C, G, E]$, making six altogether. The term ``first inversion'' is often applied loosely both to $[E, G, C]$ and $[E, C, G]$, and the term ``second inversion'' is applied to $[G, C, E]$ and $[G, E, C]$.

In the case of tetrachords there are twenty-four permutations of the way in which the four notes can be distributed over four voices. 
The minor sixth $Cm^{6} = [C, E_\flat, G, A]$, in root position, admits, among its twelve odd permutations, the form $[A, E_\flat, G, C]$.  The significance of this permutation is apparent if we transpose the root down by a major third to give  $A_\flat m^{6} = [A_\flat, C_\flat, E_\flat, F]$ in root position and its odd permutation $[F, C_\flat, E_\flat, A_\flat]$ which, in Wagner's notation 
$[F, B, E_\flat, A_\flat]$ we recognize as the ``Tristan'' chord in the form it first appears in the opening bars of {\it Tristan und Isolde}. The permutation in this case is $[1,2,3,4] \to [4,2,3,1]$. 
The same chord, at the same pitch levels, appears prominently,  in the left hand, as an arpeggio, fortissimo, in a passage of Chopin's G Minor Ballade, Op.~23, just before the {\em pi\`u animato}, in the enharmonically equivalent representation $[E_\sharp, B, D_\sharp, G_\sharp]$, obtained by applying the same permutation to the $G_\sharp m^{6}$ chord.  
We have segued on to the consideration of minor sixths here to stress the importance of the group of permutations of $n$ objects in chord formation. We shall have more to say on the Tristan chord shortly.  Note that under the permutation $[1,2,3,4] \to [4,1,2,3]$  the minor sixth $Cm^6 = [ C, E_{\flat}, G, A] $ becomes a half-diminished seventh, $A^{\varnothing 7} = [A, C, E_\flat, G]$. Similarly,  under $[1,2,3,4] \to [2,3,4,1]$ the minor seventh $Cm^7 = [C, E_{\flat}, G, B_{\flat}]$ becomes a major sixth, $E_{\flat}M^6 = [E_{\flat}, G, B_{\flat}, C]$.  Since $Cm^7$  and $E_{\flat}M^6$ are the same modulo permutation, is the underlying triad $Cm$ or $E_{\flat}M$? Likewise, since $Cm^6$  and $A^{\varnothing 7} = [A, C, E_\flat, G]$ are the same modulo permutation, is the underlying triad $Cm$ or the diminished triad $[A, C, E_\flat]$?  Such chords can perform in either capacity -- one needs the context in which they appear to determine their function. Two historical precedents for this observation are Rameau's notion of {\em double emploi} and Gottfried Weber's {\em Mehrdeutigkeit}. 

With these conventions we turn to the construction of the tonnetz.
A major triad can be changed to a minor triad in three different ways by fixing two of the tones  and altering the third. For example, in the case of $CM^{(0)}$, we can lower the $C$ by a semitone to obtain $Em^{(2)} = [B, E, G]$.  If we lower $E$ by a semitone we obtain $Cm^{(0)} = [C, E_{\flat}, G]$.  If we raise $G$ by a tone we obtain $Am^{(1)} = [C, E, A]$. 
Similarly, a minor triad can be changed to a major triad in three different ways by fixing two of the tones and altering the third. In the case of $Cm^{(0)} = [C, E_{\flat}, G]$, we can lower the $C$ by a tone to obtain $E_{\flat}M^{(2)} = [B_{\flat}, E_{\flat}, G]$ or we can raise $E_{\flat}$ by a semitone to obtain $CM^{(0)} = [C, E, G]$ or we can raise $G$ by a semitone to obtain $A_{\flat}M^{(1)} = [C, E_{\flat}, A_{\flat}]$. These are the only major triads that can be obtained from $Cm$ by fixing two of the tones and altering the third. 

It follows that associated to each unordered major triad there are three unordered minor triads, each sharing two of the pitches of the major triad; and associated to each unordered minor triad there are three unordered major triads, each sharing two of the pitches of the minor triad.
We do not assume {\em a priori} that the admissible transformations between major and minor triads are of a tone or semitone in magnitude; this is an emergent feature that arises from our assumption that each transformation should preserve two of the three tones.
%}.

In this paper we highlight -- for the first time, we believe -- the unifying role of Levi graphs and configurations in the theory of music.  
In Section II we show that the tonnetz can be represented by a finite bipartite graph of degree three and girth six \footnote{We follow the approach of Waller (1978). This fine piece of work seems to have been neglected and we have not seen Waller's bipartite graph discussed much elsewhere in the literature (but see Cohn 2012 at page 18). There is no mention of the tonnetz in Waller's paper, for at that point the tonnetz had not become an object of renewed interest to music theorists, which began with the work of Lewin and his school in the 1980s. Nor does Waller mention that he is working with a Levi graph. This is consistent with the fact that, with a few exceptions, such as Hilbert and Cohn-Vossen (1952) and Coxeter (1950), the theory of configurations, which had been active in the 19th century, had been dormant from 1910 to 1990 (Gr\"unbaum 2009, pages 8-14). Derek Waller, a British mathematician who worked in category theory, algebra and graph theory, died in 1978 at the age of 37. The connection between tonnetze and configurations is our idea, as is the realization that Waller's construction, which we now recognize as the finite form of the dual of the Eulerian tonnetz, is the Levi graph of the D222 of Daublebsky von Sterneck.}. This result is stated in Proposition 1 and the Levi graph is drawn in Figure \ref{eulerian_tonnetz}. The white vertices represent major triads; the black vertices represent minor triads. 
We show how properties of the four hexatonic cycles and the three octatonic cycles can be understood in the setting of this Levi graph and we comment on the role of other cycles, which we illustrate in Figure \ref{fig:cycles}  with several examples.
In Section \ref{sec:The Tonnetz as a Configuration} we show how the tonnetz can be regarded as a configuration. This is achieved by treating the transformations between major and minor triads as incidence relations in a combinatorial geometry.
We reach the surprising conclusion that the tonnetz can be realized geometrically in $\mathbb R^2$ as a certain self-dual configuration of lines and points. In particular, the tonnetz turns out to be the configuration $\{12_3\}$ of Daublebsky von Sterneck type D222. The result is stated in Proposition  \ref{tonnetz configuration} and the configuration is drawn in Figure \ref{fig:D222}. 

A central argument that we maintain throughout this discussion is that Euler's tonnetz is a mathematical object --  specifically, a combinatorial configuration that can be represented both by a Levi graph and a figure of lines and points in the Euclidean plane. But it also admits a representation in terms of audible relations in music. Thus, a certain abstract structure manifests itself  (i) as a combinatorial configuration, (ii) as a geometric configuration, and also (iii) as a ``musical'' configuration, the latter being the tonnetz. 

As an application of our analysis of the various cycles of the Eulerian tonnetz, in Section \ref{Cyclic Structure as a Basis for Musical Analysis} we take a look  at the problem of the resolution of the Tristan chord as it arises in  works by Wagner, Chopin and Tchaikovsky, excerpts of which we exhibit in Figures \ref{fig:Tristan}, \ref{Chopin Ballad}, \ref{tchaikovsky} and \ref{parsifal entrance}. 
In this section our analysis of compositions by these nineteenth-century composers involves a reduction from tetradic harmony to triadic harmony. Such a reduction, which can be regarded as a kind of reconnaissance, a survey of the landscape before entering into tetradic territory, is surprisingly effective and various insights can be gained on the basis of this approach. 
The analysis of tetradic harmonies by use of triadic relations can be found in the works of many authors -- and to that extent we but follow in their footsteps -- nonetheless, {\em the idea of using the Levi graph associated with the tonnetz configuration as a fundamental  basis for such analysis is entirely new}, as best as we can say, and there are three immediate benefits: (a) the discussion can often be simplified and/or made more transparent by use of the Levi graph as a highly effective ``road map'' of the tonnetz, (b) the structure of the chord space can be quite unambiguously visualized and surveyed, and (c) the approach can be generalized in a satisfying way to tetradic harmonies, which we discuss in Section \ref{sec:Tone Networks for Tetrachords}. 

At first thought it may seem puzzling that we should choose to analyze the music of Wagner and other composers of that period using both the Eulerian tonnetz, which is based on triadic harmonies, as well as the tetradic tonnetz that we develop later. Is it not inconsistent to be using two different such systems of analysis for the same purpose? Our reply is that one needs to consider the idea that Wagner and his colleagues were themselves working, in effect, with two or perhaps even three different harmonic systems in their compositional methods: first, the system of classical harmony, of which of course they were masters; but also (perhaps unconsciously) the system of the tonnetz, which somehow can be made to operate in parallel with the classical system (the literature of the tonnetz provides an abundance of examples of what we now call pan-triadic harmony); and finally, a second tonnetz system, specifically adapted to the tetradic harmonies typical of the later  works of Wagner and composers who were under his influence, concerning which we shall say more shortly. 

 Our strategy for the development of a geometric tonnetz for tetrachords is first to look in detail at the relation between major and minor triads when exactly two tones are allowed to change in a chord transition. This is carried out in Section \ref{sec:Tone Networks and Tessellations}.  We show that the resulting bipartite graph has two components, each containing six major triads and six minor triads. These graphs are shown in Figure \ref{Archimedean tonnetz}. Each component has girth four: hence neither is a Levi graph -- but each component can be unfolded to yield the infinite graph of an Archimedean or {\em semi-regular} tessellation -- a tiling of the plane by hexagons, squares, and dodecagons,  shown in Figure \ref{dodecagontessellation}. The resulting cycles lead to interesting sonorities and in Figure \ref{fanfare} we present a short composition in the form of a fanfare based on a hexacycle of this tonnetz.  

In Section \ref{sec:Tone Networks for Tetrachords} we turn our attention to the Tristan genus of tetrachords comprising the twelve dominant sevenths and the twelve minor sixths. For the term ``Tristan genus'' see Cohn (2012), pages 159-166. Our use of the term ``minor sixth'' for the sonority of the Tristan chord $[F, B, E_\flat, A_\flat]$ deserves comment. Of the 24 permutations of its constituents, two are particularly relevant. The first is the minor sixth ${A_\flat}m^6 = [A_\flat, C_\flat, E_\flat, F]$ or equivalently ${G_\sharp}m^6 = [G_\sharp, B, D_\sharp, E_\sharp] $. The second is the half-diminished seventh $F^{\varnothing 7} = [F, A_\flat, C_\flat, E_\flat]$  or equivalently ${E_\sharp}^{\varnothing 7} = [E_\sharp, G_\sharp, B, D_\sharp]$. By a minor sixth chord we mean a minor triad with an added sixth, and by a half-diminished seventh chord we mean a minor seventh chord with a flatted fifth. To suggest (whether coming consciously or otherwise from what Cohn 2012 at page 143 calls the ``third-stacking monolithy'') that the half-diminished seventh description is somehow more correct is problematic. For in {\it Tristan und Isolde} the sense of this sonority is very often explicitly that of a minor triad with an added sixth (see Bailey 1985 at pages 122-125 and Cohn 2012 at pages 142-145). Hence, we mainly use the language of minor sixths. To translate, note that the root of a minor sixth chord is a minor third higher than the root of the corresponding half-diminished seventh, so ${A_\flat}m^6  \approx F^{\varnothing 7}$.  
Previous work in this area (Childs 1998, Douthett and Steinbach 1998, Tymoczko 2011, Cohn 2012) has often tended to look at irregular graphs, with the dominant sevenths and minor sixths as vertices, and the edge sets including all possible ``parsimonious'' transformations, or some large subset thereof, variously supplemented by further vertices representing additional  tetrachords, some less consonant, such as the twelve minor sevenths, the six French sixths, and the three diminished sevenths. 

The irregularities of these graphs, the large numbers of edges, and the heterogeneity of the vertex constituents point to their inadequacy as mathematical models. Can one do better? 
As a way forward, following the combinatorial idea of the Archimedean tonnetz constructed in the previous section we consider the situation when  two tones of a tetrachord are allowed to change in a chord transition between two chords of the reductively opposite mode when we regard the minor sixth chord as a minor triad with an added dissonant sixth tone and we regard the dominant seventh chord as a consonant major triad with an added flat seventh tone. 
This gives a six-to-six map between the two sets of twelve chords. Such a map cannot give rise to a configuration, since the fundamental inequality is violated. 

To give a bipartite graph of degree six with no cycles of length less than six, the number of vertices must be no less than thirty-one -- but we only have twenty-four.  Thus, to obtain a tractable subclass of transformations that leads to a geometry, we must cut down the number of allowed transformations. Our proposed solution to the problem, which is new as far as we can see, is as follows: we constrain the transformation class in such a way that it always retains the dissonant tone and hence always changes two tones of the underlying major or minor triad. Note that the image of the dissonant tone need not itself (but might) be a dissonant tone of the transformed chord. 
The result is a three-to-three map between the dominant sevenths and the minor sixths with no tetracycles, that is to say, a configuration $\{12_3\}$ -- and we draw its Levi graph in Figure \ref{fig:tristan_tonnetz}. 

One can tell at once that this Levi graph is distinct from that of the Eulerian tonnetz by counting the numbers of cycles of various lengths. Any discrepancy suffices to show that the new tonnetz for tetrachords is different from the Eulerian one. But all the configurations $\{12_3\}$ are known, so if it is not the Eulerian tonnetz, which one is it? 

Surprisingly, it turns out that the Tristan-genus  tonnetz is the D228 of Daublebsky von Sterneck, which we depict in Figure \ref{tristan_tonnetz_configuration} -- with tetrachord labels attached.  This highly symmetric, attractive tonnetz configuration offers a new means for navigating the Tristan-genus universe. We conclude with the analysis of passages in {\it Tristan und Isolde}, {\it G\"otterd\"ammerung} and {\it Parsifal}, which allow for comparison of our methods with arguments of Lewin (1996), Childs (1998), Douthett and Steinbach (1998), Tymoczko (2011), Cohn (2012), and others.  

%%%%%%
%SECTION II
\section{The Tonnetz as a Levi Graph}
%%%%%%

\noindent The tonnetz has been studied extensively in the literature of music and mathematics. It would be out of place to attempt a bibliography here, but we mention  Lewin (1987), Hyer (1989), Cohn (1996, 1997, 1998, 2012), Douthett and Steinbach (1998), Childs (1998), Gollin (1998), Gollin and Rehding (2011), Tymoczko (2011, 2012), Catanzaro (2011),  Yust (2018), Cubarsi (2024) and Rietsch (2024) as examples of recent work  on the tonnetz.  
 
Most investigations of the tonnetz have tended to develop ideas around the triangular graph of Figure \ref{Triangular_tessellation}, representing triads as ``assemblies of atomistic pitch-class components''  (see Cohn 2012 at page  114) 
or else the ``fused-triad'' hexagonal graph of Figure \ref{tonnetz_tessellation}, along with the group-theoretic ideas implicit in both of these representations. 

Though not inconsistent with such programs, we prefer a different approach to the tonnetz -- one that allows more readily for the far-reaching generalizations ultimately required for treating other harmonic structures, such as those arising with tetrachords or with forms of music outside of the extended common practice. 

\begin{figure}[htbp] %eulerian tonnetz
\centering
\includegraphics[clip,scale=0.95]{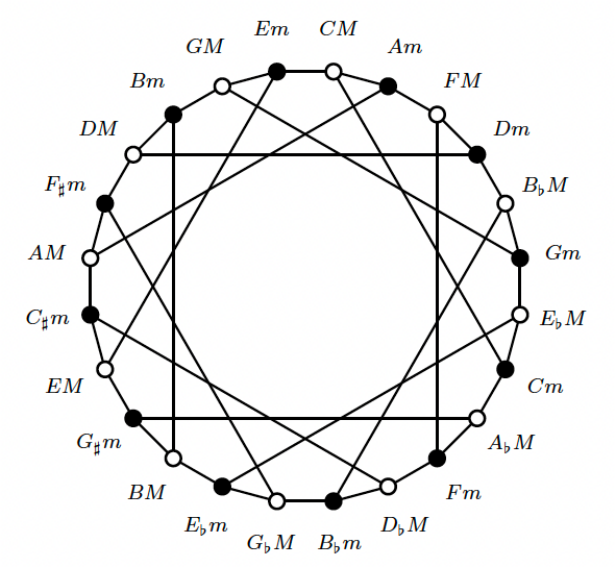}
\caption{The tonnetz as a Levi graph. The tonnetz can be modelled as a regular bipartite graph of degree three with twelve white vertices corresponding to major triads and twelve black vertices corresponding to minor triads. Each major triad is linked to three minor triads and each minor triad is linked to three major triads. The four $3p$-hexacycles are blunted triangles such as $\langle CM, Cm, A_{\flat}M, G_{\sharp}m, EM, Em, CM\rangle$. There are twelve $2p$-hexacycles (bow ties), one beginning at each major triad, an example being $\langle CM, Am, FM, Fm, A_{\flat}M, Cm, CM\rangle$. }
\label{eulerian_tonnetz}
\end{figure}

\begin{figure}[htbp] %four cycles picture
\centering
\includegraphics[clip,scale=0.70]{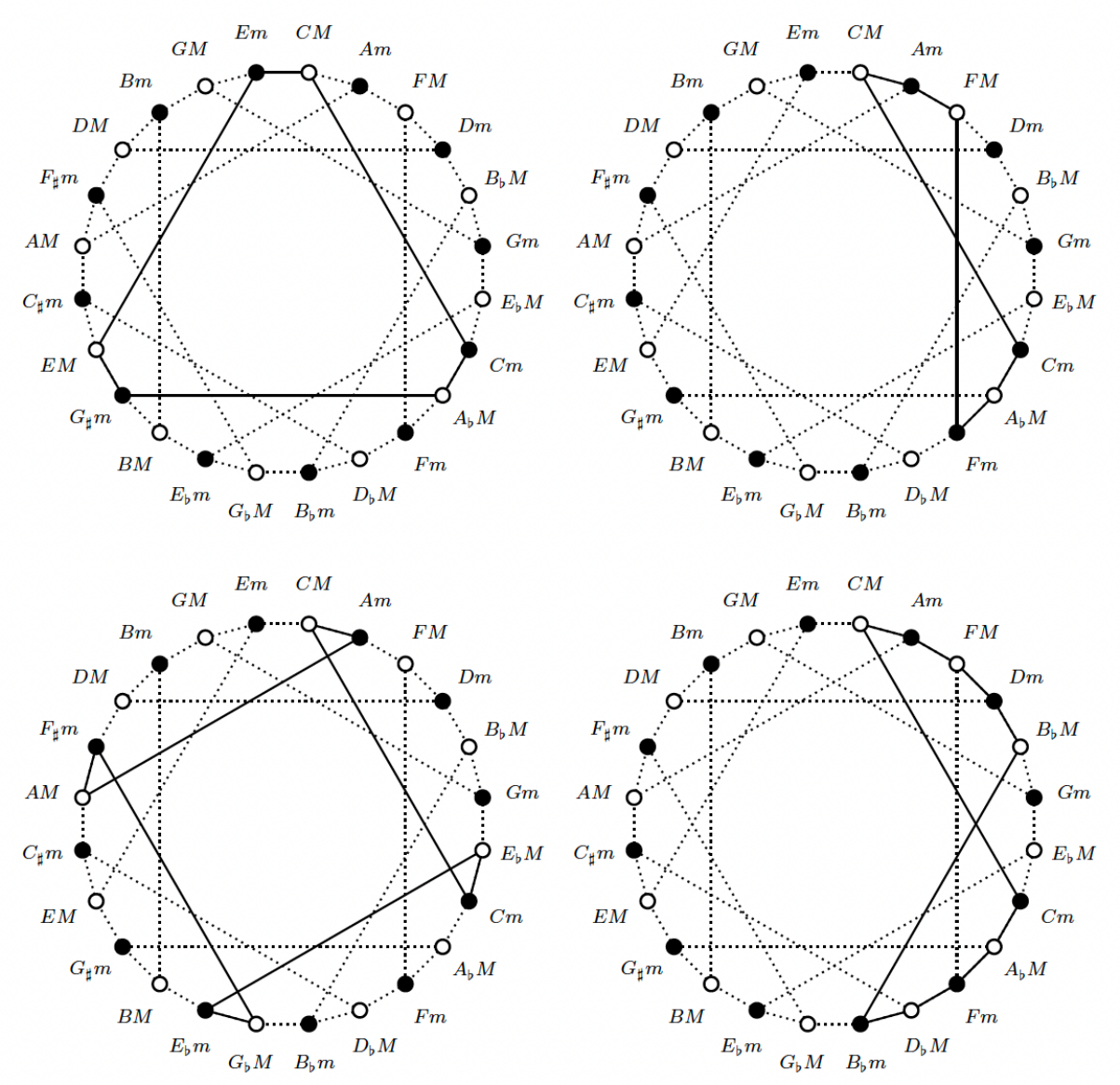}

\caption{We illustrate a $3p$-hexacycle, a $2p$-hexacycle, a $4p$-octacycle and a $2p$-decacycle. The $3p$-hexacycle $\langle CM, Cm, A_{\flat}M, G_{\sharp}m, EM, Em, CM\rangle$ is one of the four basic hexatonic sequences. The sequence  $\langle CM, Am, FM, Fm, A_{\flat}M, Cm, CM\rangle$ is one of the twelve overlapping ``straight bow tie'' $2p$-hexacycles.  Similarly,  $\langle CM, Am, AM, F_{\sharp}m, G_{\flat} M, E_{\flat} m, E_{\flat} M, Cm, CM\rangle$ is one of the three ``four-cornered hat''  $4p$-octacycles. $\langle CM, Am, FM, Dm, B_{\flat} M, B_{\flat} m, D_{\flat} M, Fm, A_{\flat} M, Cm, CM\rangle$ is one of the twelve  ``floppy bow tie''  $2p$-decacycles.  
}
\label{fig:cycles}
\end{figure}

Our starting point is the observation  that the tonnetz can be represented by the regular bipartite graph of degree three and order 24 shown in Figure \ref{eulerian_tonnetz}. Major triads are represented by white vertices and minor triads by black vertices. Each major triad is connected to the vertices of three minor triads  and each minor triad is connected to the vertices of three major triads. For example, $CM$ is connected to $Em$, $Cm$, $Am$, whereas $Cm$ is connected to $E_{\flat}M$, $CM$, $A_{\flat}M$. 
Various well-known features of the tonnetz can be read off from this graph and less well-known aspects of the tonnetz also become apparent.  

To develop our approach further, we recall a few ideas from graph theory  (Harary 1969, Wilson 1972, Bondy and Murty 1976). 
We write 
$\mathbb N = \{1, 2, 3, . \, .\, . \,\}$ for the natural numbers and $\mathbb Z^+ = \{0, 1, 2, . \, .\, . \,\}$  for the non-negative integers. There are, of course, diverse notations for these basic sets; here we follow Welsh (1976).  A graph is a triple $G = (V, E, \pi)$ where $V$ is a non-empty finite set of elements called vertices, $E$ is a set of elements called edges and $\pi: E \to \mathcal P(V)$ is a map from $E$ to the set $\mathcal P(V) = \{\{u, v\} :  u, v \in V\}$. 

We observe that $\mathcal P(V)$ is the set of unordered pairs of elements of $V$.  Care is required when we consider elements of $\mathcal P(V)$ of the form $\{v, v\}$, $v \in V$. Since $\{v, v\} = \{v\}$ as sets, we have $\mathcal P(V) = \{\{u, v\} :  u, v \in V, u\neq v \} \cup \{ \{v\} :  v \in V\}$. 
If  $\pi(e) = \{v_1,v_2\}$, where $v_1, v_2 \in V$ and $e \in E$,  we say that $v_1$ and $v_2$ are joined by $e$ and  that $e$ is incident with $v_1$ and $v_2$.  
%%%%%%%
%%%%%%
If $e_1,e_2 \in E$ are such that $\pi(e_1) = \{v_1,v_2\}$ and $\pi(e_2) = \{v_1,v_2\}$, where $v_1,v_2 \in V$, then we say that $e_1$ and $e_2$ are multiple edges. A graph $G = (V, E, \pi)$ has no multiple edges if and only if $\pi$ is injective. 
If $e \in E$ is an edge such that $\pi(e) = \{v,v\}$ for some $v \in V$, then $e$ is a loop at $v$.  If $e_1,e_2 \in E$,  $e_1 \neq e_2$, are such that $\pi(e_1) = \{v,v\}$ and $\pi(e_2) = \{v,v\}$ for some $v \in V$, then we say that $e_1$ and $e_2$ are multiple loops at $v$.

The number of edges joining a vertex is the {\it degree} of that vertex. The degree of a graph is the maximum of the degrees of its vertices. If all the vertices have the same degree,  the graph is said to be  {\it regular}. 
If the edges $e_1$ and $e_2$ are incident with a common vertex, we say that $e_1$ and $e_2$ are adjacent. If a vertex is incident with no edges it is isolated. If $G$ has no loops or multiple edges it is said to be simple. 
A sequence of $k \in \mathbb N$ distinct consecutively adjacent edges is a path of  length $k$. A graph is  connected if any two vertices are joined by some path. 

A subgraph of $G = (V, E, \pi)$ is a graph $G' = (V', E', \pi')$ such that $V' \subset V$, $E' \subset E$, $E' = \pi^{-1}(\mathcal P(V'))$, and $\pi'$ is the restriction of $\pi$ to $E'$.  If $G$ is not connected then it splits into two or more connected subgraphs called components of $G$. 
A closed sequence of $k$ distinct consecutively adjacent edges that begins at some vertex $v$ and comes back to $v$ in such a way that all $k$ of the vertices are distinct is called a $k$-cycle. A loop is a $1$-cycle; a pair of multiple edges form a $2$-cycle.
The least value of $k$ for which a graph $G$ admits a $k$-cycle  is called its girth.  If $V$ has $N$ elements then $N$ is the order of the graph. If the vertices of a regular graph are of degree $d$, then $dN = 2M$, where $M$ is the number of edges. By a  forest we mean a graph with no cycles. By a tree we mean a connected forest. 
Finally, we say that a graph is  bipartite if sets $V_0(G)$ and $V_1(G)$ exist satisfying $V(G) = V_0(G) \cup V_1(G)$ and $ V_0(G) \cap V_1(G) = \varnothing$ such that every edge of $G$ joins a vertex of  $V_0$ to a vertex of  $V_1$.  A graph is bipartite if and only if every cycle has an even number of edges.

We also require the combinatorial idea of an incidence structure. By an incidence structure we mean a disjoint pair of finite sets $\mathbb P$ and $\mathbb L$, whose elements we call  points and lines, satisfying certain conditions. We assume that $\mathbb P$ and $\mathbb L$ are equipped with an incidence relation 
$R: \mathbb P \times \mathbb L \to \{0,1\}$
such that for any $p \in \mathbb P$ and any $L \in \mathbb L$ either $R(p,L) = 1$, when we say that $p$ and $L$ are incident, or $R(p,L) = 0$, when they are not incident. To aid intuition, we use geometrical terms such as ``$p$ lies on $L$'' or ``the lines $L$ and $M$ intersect at the point $p$'' to express incidence relations. The use of such language  is purely for convenience; since the terms ``line'' and ``point'' have no {\em a priori} geometrical meaning. We assume that two distinct lines intersect at most once and that two distinct points lie on at most one line in common; thus, if $L, M \in \mathbb L$, $L \neq M$, then $\exists$ at most one $p \in \mathbb P$ such that $R(p,L) = 1$ and $R(p, M) = 1$; if $p, q \in \mathbb P$, $p \neq q$, then $\exists$ at most one $L \in \mathbb L$ such that $R(p,L) = 1$ and $R(q, L) = 1$.
We call a triple $\{\mathbb P, \mathbb L, R\}$ satisfying these conditions an incidence structure. 

By a combinatorial configuration of type $\{ m_r, n_k \}$ for $m, n, r, k \in \mathbb N$ we mean an incidence structure with $m$ points and $n$ lines such that $k$ points lie on each line and $r$ lines go through each point. 
We use the term ``type'' here since in some cases several different configurations of the same type may exist.  
By a flag we mean a pair $\{p, L\}$ where $p \in \mathbb P$ and $L \in \mathbb L$. 
For any configuration $\{ m_r, n_k \}$ we have the fundamental equality
$mr = nk$,
since each side of this equation gives a way of counting the total number of distinct flags in the configuration.  The fundamental inequality 
\begin{eqnarray}
m \geq r(k-1) + 1
\label{fundamental inequality}
\end{eqnarray}
follows from the fact that the total number of points $m$ must be no less than the total number of points lying on any given set of flags sharing a point in common.  
If the number of lines equals the number of points in a configuration then by the fundamental equality the number of points on each line must equal the number of lines through each point, and we see that the configuration is balanced.  

There is a special type of graph associated with any combinatorial configuration called the ``Levi graph'' of that configuration (Levi 1929, 1942, Coxeter 1950). The Levi graph $G(C)$ of a combinatorial configuration $\mathscr C$ of type $\{m_r,n_k\}$ is the bipartite graph that has $m$ white vertices, corresponding to points of $\mathscr C$, and $n$ black vertices, corresponding to lines of $\mathscr C$. 
A pair of vertices of $G$ determine an edge of the graph if and only if one represents a point of $\mathscr C$ and the other represents a line of $\mathscr C$ that is incident with that point. It follows that (i) each edge  joins a white vertex and a black vertex, (ii) $r$ edges meet each white vertex, and (iii) $k$ edges meet each black vertex. Clearly, for such a graph the total number of edges is given by $mr = nk$.

The significance of Levi graphs is that their correspondence with combinatorial configurations is one-to-one (Gr\"unbaum 2009, page  28). Specifically,  a graph $G$ is the Levi graph of a  combinatorial configuration $\mathscr C$ of type $\{ m_r, n_k \}$  if and only  if (i) $G$ is bipartite with partition $V(G) = V_0(G) \cup V_1(G)$, where $V_0$ has cardinality $m$ and $ V_1$ has cardinality $n$; (ii) every vertex in $V_0$ is of degree $r$ and every vertex in $V_1$ is of degree $k$;  (iii) $G$ is simple;  (iv) $G$ has girth at least six. 
With these definitions and facts at hand, we are able to make the following assertion. 

\begin{Proposition}
The tonnetz can be represented by the Levi graph of a configuration of type $\{12_3\}$  in which the twelve major triads are represented by white vertices and the twelve minor triads are represented by black vertices. The $36$ edges are arranged in such a way that each white vertex is connected to three black vertices and each black vertex is connected to three white vertices -- that is to say, the graph is regular of degree three. 
\end{Proposition}

Inspection of Figure \ref{eulerian_tonnetz} shows that the minimal cycles are hexacycles. Hexacycles can be classified by the number of parallel transformations that they contain.  Figure \ref{eulerian_tonnetz} shows that there are exactly four hexacycles with three parallel transformations, each taking the form of a blunted triangle. 
These hexacycles partition the tonnetz and are given by
\begin{eqnarray}
\hspace{0.30cm} \langle CM, Cm, A_{\flat}M, G_{\sharp}m, EM, Em, CM \rangle, \hspace{0.5cm} 
\langle FM, Fm, D_{\flat}M, C_{\sharp}m, AM, Am, FM \rangle, \hspace{0.5cm} \nonumber
\end{eqnarray}
\vspace{-1cm}
\begin{eqnarray}
\langle B_{\flat}M, B_{\flat}m, G_{\flat}M, F_{\sharp}m, DM, Dm, B_{\flat}M \rangle, \hspace{0.5cm} 
\langle E_{\flat}M, E_{\flat}m, BM, Bm, GM, Gm, E_{\flat}M \rangle.
\end{eqnarray}
Recall that when we write out a sequence of triads forming a cycle, we use angled brackets and we repeat the initial triad at the end of the sequence to show how it closes up.  One recognizes these sequences  as the four mutually disjoint  {\it hexatonic cycles} arising in pan-triadic theories (Cohn 1996, 1997, 1998, 2012, Douthett and Steinbach 1998). See Cohn (2012) at page  211 for a definition of ``hexatonic cycle''.

Each such hexacycle is constructed from  six distinct tones and contains three distinct parallel transformations represented by moving across the Levi graph. 
We call these cycles {\it thrice-parallel}, which we abbreviate as $3p$. One of the advantages of the use of a Levi graph for the tonnetz is that it allows one to see at a glance the symmetries of the $3p$-hexacycles and how they partition the tonnetz. 

The edges of the tonnetz fall into three sets of twelve. In the first set we have edges of the form $\{CM, Cm\}$, which we already mentioned, that transform a major ({\em resp.}, minor) triad into its parallel minor ({\em resp.}, major), the $\rm\bf P$ transformations. In the second set we have edges of the form $\{CM, Am\}$ that transform a major ({\em resp.}, minor) triad into its relative minor ({\em resp.}, major), the $\rm\bf R$ transformations. In the third set  we have edges of the form $\{CM, Em\}$ that transform a major ({\em resp.}, minor) triad into a corresponding minor ({\em resp.}, major) under a leading tone exchange, the $\rm\bf L$ transformations. Three edges meet at each vertex, one of each type. 

Apart from the four hexatonic cycles, there are twelve further hexacycles in the Levi graph of the tonnetz, one beginning at each major triad, of which a typical example is
\begin{eqnarray} \label{CM twice parallel cycle}
\langle CM, Am, FM, Fm, A_{\flat}M, Cm, CM \rangle.
\end{eqnarray}
This cycle includes only  two parallel transformations -- viz., from $CM$ to $Cm$ and from $FM$ to $Fm$; hence, the twelve $2p$-hexacycles are distinct from the four $3p$-hexacycles. 
An $\rm\bf L$ transformation sends a major triad one notch counterclockwise on the Levi graph to the previous minor triad, and it sends a minor triad one step clockwise to the next major triad.  An $\rm\bf R$ transformation sends a major triad a notch clockwise on the Levi graph to the next minor triad, and it sends a minor triad one step counterclockwise to the previous major triad.  A $\rm\bf P$ transformation sends any triad across the Levi graph to its mode reverse.
The $2p$-hexacycles are what Cohn (2012), pages 113-121, calls ``pitch retention loops''; at page  43 of Cohn (1997), they are called ``$\rm\bf LPR$ loops''.

%Section III
\section{The Tonnetz as a Configuration}
\label{sec:The Tonnetz as a Configuration}
%%%%%%

\begin{quote}  \hspace{1cm} ``.\,.\,.\,interesting configurations are represented by interesting graphs.''

 \hfill ---H.~S.~M.~Coxeter,~{\em Twelve Geometrical Essays}
\end{quote}
\vspace{0.15cm}

\noindent If we accept the view that Euler's tonnetz is not merely a curiosity but is fundamental to the structure of music of the common practice period, then it makes sense to look at the mathematics of the tonnetz from every angle -- group-theoretic, number-theoretic, set-theoretic, combinatoric, and geometric -- for  insights. 
One key concept, at a nexus of the various branches of mathematics mentioned above, is the idea of configuration that we introduced in Section \ref{sec:Introduction}. 
By a configuration in the broadest sense, by which we mean combinatorial, we mean a 
set of $m$ points and $n$ lines such that every point lies on $r$ lines and every line contains $k$ points. Variants of the idea arise when we consider the interpretations of the words ``point'' and ``line'' and when we consider the spaces in which these elements are taken. 

Configurations $\{n_3\}$ have attracted much attention and have been investigated at length. The inequality \eqref{fundamental inequality} shows that no $\{n_3\}$ configurations exist for $n \leq 6$.  We draw on this fact a little later.
Let us take stock of some of the configurations that exist for higher values of $n$. 
There is a unique configuration $\{7_3\}$,  modelled by the seven points and seven lines of the finite projective plane based on the binary field $\mathbb Z_2$. This is the Fano plane.
There is a unique combinatorial configuration $\{8_3\}$, which can be realized explicitly by a system of eight points and eight complex projective lines in the complex projective plane. No real coordinatization of $\{8_3\}$ on the Euclidean plane exists. 
There are three distinct combinatorial configurations $\{9_3\}$, all of which can be modelled as systems of points and lines in $\mathbb {R}^2$.  One of these is the well-known Pappus configuration. 
 All three of the configurations $\{9_3\}$ can be modelled over the rational numbers. 
 
 There are ten distinct  combinatorial configurations $\{10_3\}$, nine of which can be modelled as geometric configurations in 
 $\mathbb {R}^2$ and are rationally realizable. Among these is the famous Desargues configuration, which can be used to construct a tonnetz for pentatonic music.  But one of the combinatorial configurations  $\{10_3\}$  does not admit a representation as a geometric configuration over any field.  
The 31 distinct configurations $\{11_3\}$ were identified by Daublebsky von Sterneck (1894), all of which can be realized in  $\mathbb {R}^2$.  

The situation with $\{12_3\}$ is not completely straightforward. 
Daublebsky von Sterneck (1895) established the existence of 228 distinct configurations of this type. This number was accepted even as recently as the work of Sturmfels and White (1990), who were investigating Gr\"unbaum's conjecture that every $\{n_3\}$ that can be realized over the reals can be realized over the rationals. 
They established by use of diophantine computations that all the 
$\{11_3\}$ and $\{12_3\}$ cases considered by Daublebsky von Sterneck were realizable over $\mathbb Q$.
Then Gropp (1990, 1993, \!1997) showed that one case was missing from Daublebsky von Sterneck's list:  in fact, the number of combinatorial configurations $\{12_3\}$  is 229. 
Thus it follows from a long line of work by geometers and combinatorialists spanning more than a century that all 229 of the combinatorial $\{12_3\}$ configurations admit geometric realizations in  $\mathbb R^2$ and indeed in the rational plane. See Betten, Brinkmann and Pisanski (2000), Gropp (2004), Gr\"unbaum (2009), and Alazemi and Betten (2014) for the history of the subject and for tables and constructions for $\{12_3 \}$ and $\{ n_3 \}$ for higher values of $n$.
A configuration is said to be ``self-dual''  if when we interchange the roles of lines and points we get a configuration with the same incidence relations as the original. 
The triangle $\{3_2\}$ is self-dual and for each $m \geq 3$ we can construct a self-dual configuration $\{m_2\}$ in the plane. It should be evident that any self-dual configuration is balanced. On the other hand, a balanced configuration is not necessarily self-dual. This is not so obvious, since the combinatorial configurations associated with low order symmetric Levi graphs are self-dual. 
Among the balanced combinatorial configurations, the Fano $\{7_3\}$, the M\"obius
$\{8_3\}$, the Pappus $\{9_3\}$,  the Desargues $\{10_3\}$, along with the other two $\{9_3\}$ configurations and the nine other $\{10_3\}$ configurations are self-dual. But of the 31 configurations of type $\{11_3\}$, only 25 are self-dual, and of the 229 configurations of type $\{12_3\}$, only 95 are self-dual (Betten, Brinkmann and Pisanski 2000). 

\begin{figure}[htbp] %D222
\centering
\includegraphics[clip,scale=0.60]{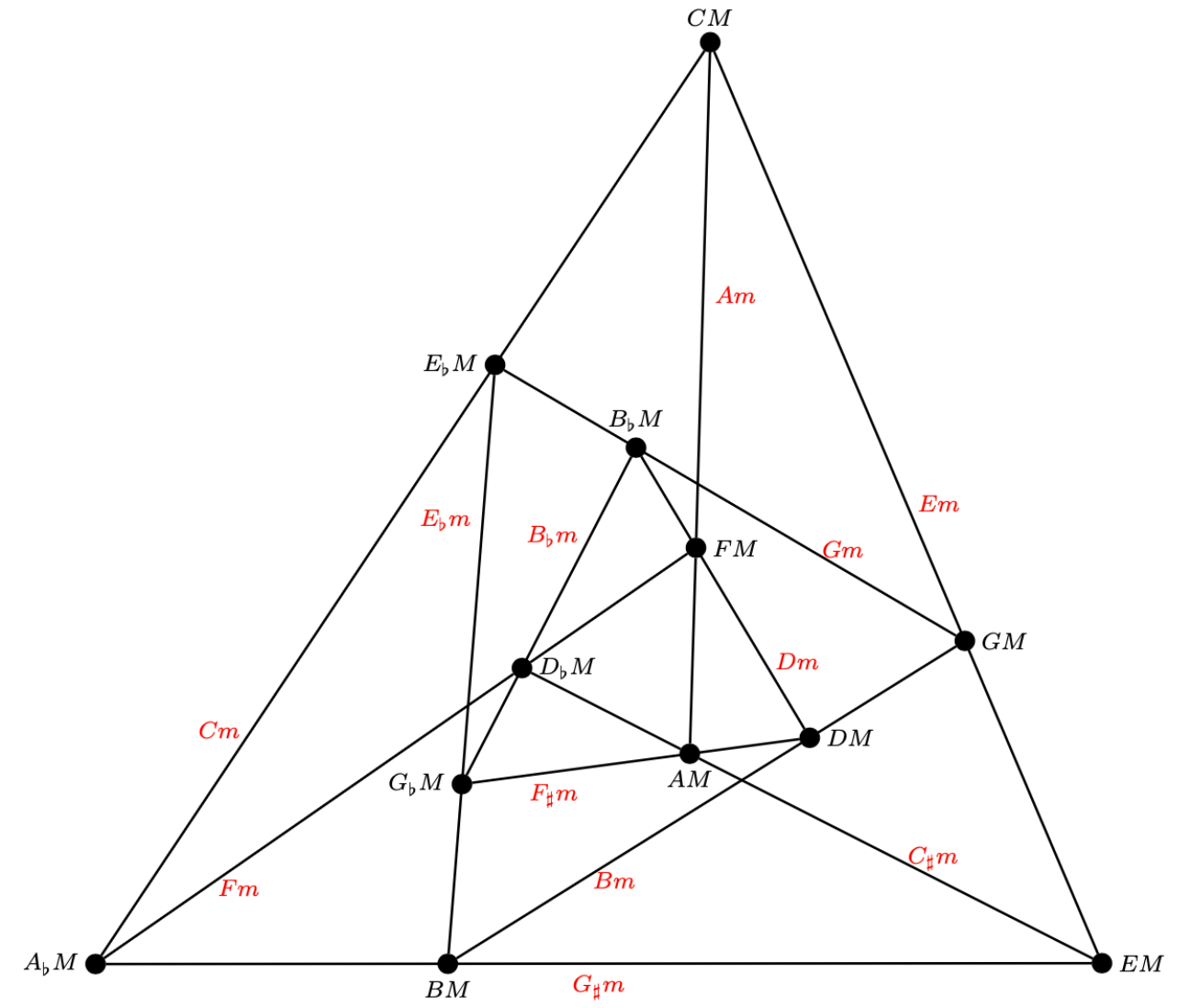}

  \caption{The tonnetz admits a realization as a configuration in $\mathbb R^2$ with a $\{12_3\}$ incidence scheme of Daublebsky von Sterneck type D222. Major triads are represented by points and minor triads by lines. Through each point representing a major triad one finds three lines representing the minor triads associated with it; and on each such line representing a minor triad one finds three points representing the three associated major triads. The four hexatonic $3p$-hexacycles such as $\langle CM, Cm, A_{\flat}M, G_{\sharp}m, EM, Em, CM \rangle$ form a closed sequence of four triangles with the property that the vertices of each such triangle are inscribed in the lines of the next.}

\label{fig:D222}
\end{figure}

In his paper on the configurations $\{11_3\}$,  Daublebsky von Sterneck (1894) makes the following remark, which we can presume applies also to his work on configurations $\{12_3\}$: ``Es handelt sich hier nur darum, alle m\"oglichen Schemata, welche obige Eigenschaften besitzen, herzustellen und ist die
Aufgabe insoferne ein Problem der Combinationslehre. Die im
Folgenden verwendeten geometrischen Betrachtungen dienen bloss
dazu, die combinatorischen \"Uberlegungen zu unterst\"utzen und zu erleichtern.'' 

In English: ``Here the only goal is to create all possible systems which possess the above $[\,\{11_3\}\,]$ properties and  in this respect the task is a problem of combinatorics. The geometric considerations used in the following serve merely to support and clarify the combinatorial considerations.''  

The remark is interesting in the present context, since it suggests that one can take a minimalist view, that the tonnetz is a combinatoric object, represented by its Levi graph, and that its geometric representation as a configuration mainly serves to support and clarify the combinatorial considerations.  

But combinatorics and geometry have come a long way since the days of Daublebsky von Sterneck. The modern perspective is to take a more unified approach to these disciplines, and we can take the view that music is deeply entangled with this unity. The proof of concept for such an outlook is simply to see whether entertaining this view is fruitful in the new ideas that it generates in the context of mathematical music theory -- and that is our goal in the present work. 

Now we are in a position to say in more detail how the tonnetz fits into this picture. We know from Figure \ref{eulerian_tonnetz} that the tonnetz can be represented as the Levi graph of a combinatorial configuration of type $\{12_3\}$. 
But we also know that each such combinatorial configuration $\{12_3\}$ can be realized as a geometric configuration 
$\{12_3\}$.  
That implies that the tonnetz is such a configuration -- that it can be drawn as a geometric configuration of points and lines in the Euclidean plane, the major triads being represented by points and the minor triads by lines, so that the points lie on lines in threes and the lines go through points in threes. 
But which configuration? Once one gets this far, the problem that remains is to determine which one of 229 known configurations of type $\{12_3\}$ is the Eulerian tonnetz. 

We have been able to identify it as D222. This remarkable configuration of points and lines, labelled with major and minor triads, which we call the tonnetz configuration, is depicted in Figure \ref{fig:D222}, in line with the original drawing of Daublebsky von Sterneck (1895). Thus, we have obtained a truly interesting and unexpected link between music and geometry: 

\begin{Proposition}
The tonnetz can be represented in  $\mathbb {R}^{2}$ as a self-dual configuration $\{12_3\}$ of Daublebsky von Sterneck type  {\em D222}. The twelve major triads are represented by points and the twelve minor triads are represented by lines. Three of the points lie on each line and three of the lines pass through each point. The thirty-six incidence relations between the twelve points and the twelve lines determine the edges of the associated Levi graph. 
\label{tonnetz configuration}
\end{Proposition}

%%%%%%
%Section IV
\section{Cyclic Structure as a Basis for Musical Analysis}
\label{Cyclic Structure as a Basis for Musical Analysis}
%%%%%%
%%%%%%
\noindent Clearly, any geometric configuration gives rise to a combinatorial configuration.
It follows that any geometric configuration gives rise to a Levi graph, the white vertices corresponding to points and the black vertices to lines. Thus a triangle in $\mathbb {R}^2$ determines a hexagonal Levi graph. 
In particular, \textit{the hexacycles of the tonnetz are represented by triangular subconfigurations of the tonnetz configuration}. The 3$p$-hexacycle containing $CM$, for example, can be represented by a triangle in which the three major triads $CM$, $EM$, $A_{\flat}M$ correspond to points and the three minor triads $Em$, $G_{\sharp}m$, $Cm$ correspond to the lines that join them. 

But the 2$p$-hexacycles also correspond to triangles in the tonnetz configuration. For example, the 2$p$-hexacycle $\langle CM, Am, FM, Fm, A_{\flat}M, Cm, CM \rangle$ is the triangle with vertices $CM, FM, A_{\flat}M$. It is a simple exercise to count the total number of triangles in Figure \ref{fig:D222}, allowing us to determine altogether sixteen different hexacycles in the tonnetz. Thus, each representation of the tonnetz -- whether it be tessellation, Levi graph, or configuration -- has merit in revealing aspects of musical structure.

 \begin{figure}[htbp] %Schubert
\centering
\includegraphics[clip,scale=0.45]{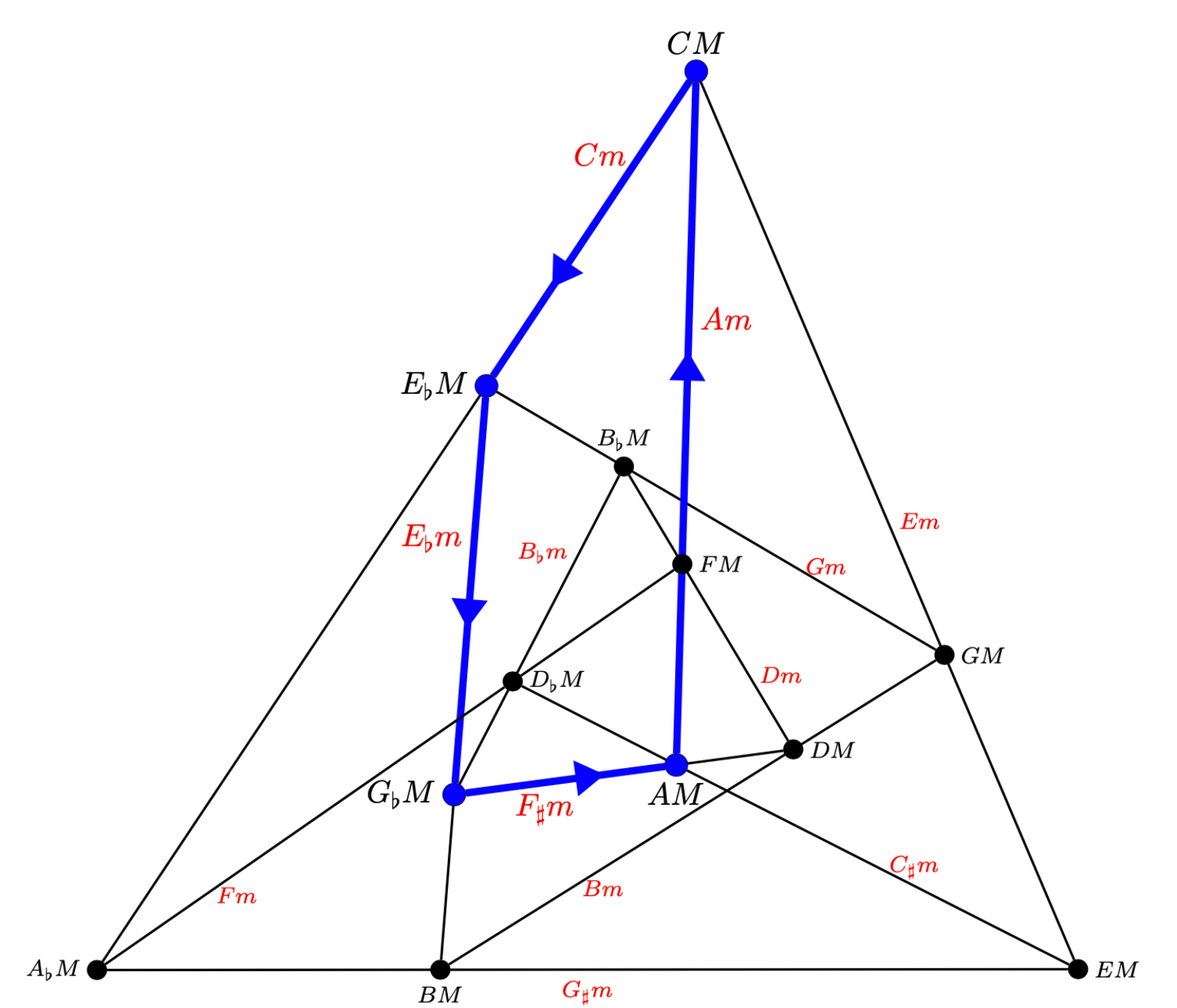}
 \caption{The 4$p$-octacycle in Schubert's Overture to {\em Die Zauberharfe}. The progression begins with the line representing $Cm$ and continues  counterclockwise to $E_{\flat}M$ and eventually on to $CM$.}
\label{fig:Schubert}
\end{figure}

Similarly, one sees that \textit {the octacycles of the tonnetz are represented by quadrilateral subconfigurations of the tonnetz configuration}.  The 4$p$-octacycle $\langle Cm$, $E_{\flat}M$, 
$E_{\flat}m$, $G_{\flat}M$, $F_{\sharp}m$, $AM$, $Am$, $CM$, $Cm \rangle$ arising as the blade-like quadrilateral subconfiguration of the tonnetz configuration shown in Figure \ref{fig:Schubert} may look familiar. It is identical to the chord sequence of the  ``four-cornered hat'' starting at $Cm$ in Figure \ref{eulerian_tonnetz}. 
It is also identical to the chord sequence in Schubert's Overture to {\em Die Zauberharfe} discussed in Cohn (1997) at page   35 and Cohn (2012) at pages  85-89.  The three 4$p$-octacycles, which correspond to a system of three such mutually exclusive blade-like ``box-cutter'' quadrilaterals,  partition the tonnetz and the trajectory in Schubert's overture is confined to one of these sequences, as one sees in Figure \ref{fig:Schubert}. So here we have a vivid geometric representation of Schubert's sequence.
This example illustrates the point that although strictly speaking there is no ``musical information'' in such a geometric configuration that is not already implicit in the associated Levi graph, since these constructions are isomorphic, nonetheless, the visualization of geometric configurations makes them ideal for the representation of certain musical ideas. Another example of this phenomenon is the surprising visualization of the octacyclic progression of minor sixth chords and dominant sevenths in the Tristan-genus tonnetz presented in Section \ref{sec:Tone Networks for Tetrachords} in our discussion of the final scene of G\"otterd\"ammerung.
 
A cycle of the tonnetz represents a sequence of chords progressing via parsimonious voice leading in such a way that the sequence eventually closes upon itself.  
One can check that each triad, major and minor, belongs to exactly three $2p$-hexacycles. 
That the twelve $2p$-hexacycles overlap in this way adds an element of richness to the scheme of triadic relations that may at first glance seem formidable. On reflection, one sees that the pattern of relations among the $2p$-hexacycles adds interest in a way that contrasts remarkably with the mutual disjointedness of the four $3p$-hexacycles, hence allowing for a variety of structural modalities. 
As an example, we mention the possibility of ``surgical insertion'' of one cycle within another. Suppose we consider the 2$p$-hexacycle \eqref{CM twice parallel cycle} alongside the transposed cycle lying a fourth above, $\langle FM, Dm, B_{\flat}M, B_{\flat}m, D_{\flat}M, Fm, FM \rangle$.
Since the short sequence $[FM, Fm]$ appears as a subsequence of the cycle \eqref{CM twice parallel cycle}, we can ``snip'' the cycles at $[FM, Fm]$, and splice them together to create a 2$p$-decacycle,
\begin{eqnarray}
\langle  CM, Am, FM, Dm, B_{\flat}M, B_{\flat}m, D_{\flat}M, Fm, A_{\flat}M, Cm, CM\rangle.
\end{eqnarray}
The result is not displeasing and is immediately suggestive of a new manner of composition whereby cycles can be combined to produce new cycles.  For example, the procedure can be iterated, allowing one to splice in the sequence
$\langle B_{\flat}M, Gm, E_{\flat}M, E_{\flat}m, G_{\flat}M, B_{\flat}m, B_{\flat}M\rangle$
by a similar snip-and-insert technique, leading to the still longer cycle
\begin{eqnarray}
\langle CM, Am, FM, Dm, B_{\flat}M, Gm, E_{\flat}M, E_{\flat}m, G_{\flat}M, B_{\flat}m, D_{\flat}M, Fm, A_{\flat}M, Cm, CM\rangle.
\end{eqnarray}

We observe, finally, that beginning at $CM$, the progression moving clockwise around the perimeter of the Levi graph involves a sequence of $\bf L$ and $\bf R$ transformations, leading us from $CM$ to $Am$ to $FM$ to $Dm$ and so on. 
In this way, we march from vertex to vertex through all twenty-four major and minor chords, giving us a Hamiltonian cycle. As Tymoczko (2012) remarks,  a Hamiltonian path touches on all the vertices in a graph without passing through any of them twice, and is in that sense a generalization of the twelve-tone row. A well-known  progression  in the Scherzo from Beethoven's Ninth Symphony (at measures 143-176) pointed out by Cohn (1992) gives an interesting example of the use of the perimeter Hamiltonian,
\begin{eqnarray}
[CM, Am, FM, Dm, B_{\flat}M, Gm, E_{\flat}M, Cm, A_{\flat}M, Fm, \nonumber\\ 
D_{\flat}M, B_{\flat}m, G_{\flat}M, E_{\flat}m, BM, G_{\sharp}m, EM, C_{\sharp}m, AM].\,\,\,\,
\end{eqnarray} 
Beethoven breaks off the sequence before it completes, a precedent followed by others, such as Wagner, when tonnetz cycles are incorporated into musical passages. Since there are no fewer than 62 distinct Hamiltonian cycles in the Eulerian tonnetz (see Appendix), each encompassing the twenty-four major and minor triads, any one of which cycles Beethoven might  have chosen for this passage in his symphony, one sees that the layout of the Levi graph in Figure \ref{eulerian_tonnetz} somehow prefigures Beethoven's thinking. In fact, no other representation of the tonnetz shows this Hamiltonian sequence so vividly. 

As usual with matters of artistic creativity, there are no fixed recipes for the modes in which compositional elements should be employed; rather, merely suggestions of possibility.  
Nonetheless, we can let mathematics lead the way. On the matter of combining cycles to produce new cycles, some clarity can be added if we recall the following well-known result from graph theory. 
Let $C_1$ and $C_2$ be cycles of a graph $G$ and suppose that the two cycles share an edge $e$ in common.  Then 
there exists a cycle $C_3 \subseteq \{C_1 \cup C_2\} \backslash e$. 
That was indeed the setup, for example, with the hexacycles that shared the edge $[FM, Fm]$, hence leading to a decacycle. In more detail, taking the union of the sets of the edges in the two hexacycles and removing the common edge, we were left with a set whose elements form a decacycle.  Examples of the $3p$-hexacycles and $2p$-hexacycles associated with the tonnetz can be found in Figure \ref{fig:cycles}, where we also find examples of  $4p$-octacycles and $2p$-decacycles.  %%%%%%%%%%%%%

We consider these cycles, alongside other such examples, to be architectural structures that can be used by composers in their work. 
The $4p$-octacycles are of interest as objects of musical analysis on account of their relationship to octatonic scales. An octatonic scale is a scale with eight tones in it, as opposed to the seven tones of the conventional major and minor scales, with the property that adjacent notes within the scale differ either by a semitone or a tone, alternating the gap at each step. 
The construction of such scales is as follows. There are three diminished seventh chords, given by $X = [C, E_{\flat}, G_{\flat}, A]$,  $Y = [D_{\flat}, E, G, B_{\flat}]$, and $Z = [D, F, A_{\flat}, B]$. If we form the union of two of these sets of four notes, this can be done in three ways, leading to three distinct sets of eight notes, the three octatonic scales:  $O_{12} = X \cup Y = [C, D_{\flat}, E_{\flat}, E, G_{\flat}, G, A, B_{\flat}, C]$, $O_{23} = Y \cup Z = [D_{\flat}, D, E,  F, G, A_{\flat},  B_{\flat}, B, D_{\flat} ]$ and $O_{31} = Z \cup X = [D, E_{\flat}, F, G_{\flat}, A_{\flat},  A, B, C, D]$.
Any octatonic scale, starting at any chosen note, belongs to one of these three types, possibly after a cyclic permutation. 

Octatonic scales have been used  by many composers, including Rimsky-Korsakov, Scriabin, Ravel, Stravinsky and Messiaen.
Their  relation to the tonnetz is as follows (Douthett and Steinbach 1998). 
Four major triads and four minor triads can be constructed from the constituents of an octatonic scale. For example, from $O_{12}$ one constructs $CM$, $Am$, $AM$, $F_{\sharp}m$, $G_{\flat}M$, $E_{\flat}m$, $E_{\flat}M$, and $Cm$. These are the triads of the 4$p$-octacycle starting at $CM$ in the tonnetz (see lower left quadrant of Figure \ref{fig:cycles}). The other two 4$p$-octacycles of the tonnetz can be constructed similarly. 
The 4$p$-octacycles are tied to the three diminished seventh chords in the way that the four $3p$-hexacycles are tied to the four augmented triads. 
By taking the six possible unions of pairs of the four augmented triads we obtain six distinct hexachords. From four of these, one constructs the four $3p$-hexacycles. The remaining two hexachords form the whole-tone scales, each of which contains six notes. 
The two whole-tone scales are given by $W_{1} = [C, D, E, G_{\flat}, A_{\flat}, B_{\flat}]$ and  $W_{2} = [D_{\flat}, E_{\flat}, F, G, A, B]$.
These arrays of six equally spaced whole notes have been exploited by composers such as Bart\'ok, Berg, Debussy, Jan\'a\v cek, Liszt, Puccini, and Rimsky-Korsakov, among others.  

Examples of melodic and harmonic structures based on cycles of the tonnetz in music of the common practice period have been noted by Cohn (2012) and others  in works of Mozart, Beethoven, Schubert, Schumann, Chopin, Wagner, Liszt, Bruckner, Brahms, Mahler and Strauss. Such use is not obvious to the classically trained ear since the progressions associated with tonnetz cycles do not form part of the system of harmony  taught in everyday music theory. 
Nor would the composers listed have been conscious of their own implicit use of  the cycles of the tonnetz.  Nonetheless, if the unconscious mind will tap the resources of the mathematical structures embedded in the geometry of the tonnetz then surely it will seek to avail itself of the full richness of these resources.  Taken collectively, the cycles of the tonnetz can act as a basis both for the composition of music and for its analysis; thus, in our view the general scheme of pan-triadic musicology should involve all the 5409 cycles of the tonnetz. 

%Tristan fragment
\begin{figure} [htbp]
%\centering
\includegraphics[clip, scale=0.75]{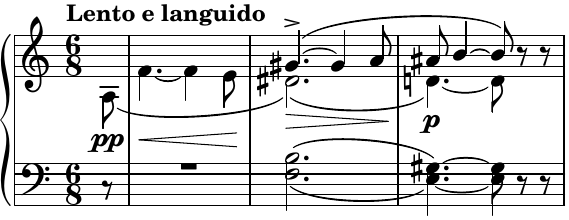}
\caption{Wagner, {\em Tristan und Isolde}, Act 1, opening bars.
%mm.\,1-3. 
The Tristan chord $[F, C_\flat, E_\flat, A_\flat]$, which in Wagner's notation takes the form $[F, B, E_\flat, A_\flat]$, is one of the twelve odd permutations of the minor sixth chord $A_\flat m^6 = [A_\flat, C_\flat, E_\flat, F]$. It can equally be regarded as an even permutation of the half-diminished seventh
$F^{\varnothing7} = [F, A_\flat, C_\flat, E_\flat]$.}
\label{fig:Tristan}
\end{figure}

%
%Chopin excerpt
\begin{figure} [htbp]
%\centering
\includegraphics[clip, scale=0.8]{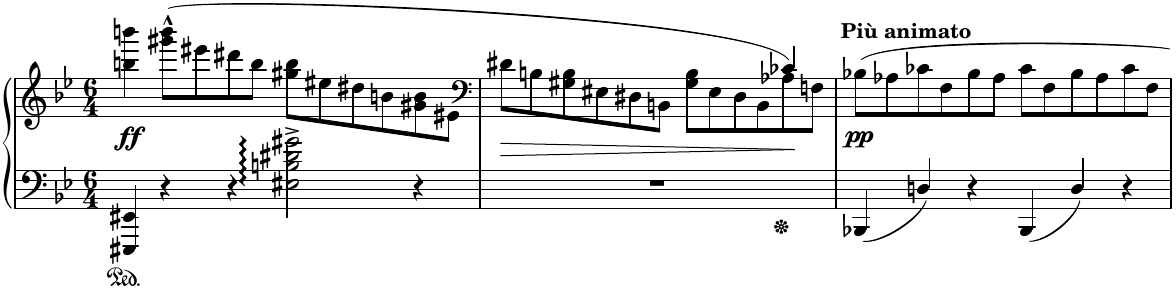}
\caption{Chopin, Ballade in G minor, Op.~23, excerpt. The Tristan chord appears  in the left hand, as an arpeggio, fortissimo, at Wagner's pitches, but in the enharmonically equivalent form $[E_\sharp, B, D_\sharp, G_\sharp]$. The descending figure in the right hand touches the same notes and is mirrored by the violins in {\it Tristan und Isolde} at the climax of the prelude. Chopin's resolution is to ${B_{\flat}}^{7}$. } 
\label{Chopin Ballad}
\end{figure}

\newpage
As an example of the analysis of cyclic structure, let us return to the resolution of the Tristan chord. Wagner offers a ``tentative'' resolution of this enharmonic $G_{\sharp}m^{6}$ to $E^{7}$, as shown in Figure \ref{fig:Tristan}. Chopin, on the other hand, resolves the same chord to $B_{\flat}^{7}$, as we see in Figure \ref{Chopin Ballad}. Chopin's resolution is definitive: it brings that lengthy first section of the Ballade to a close, and paves the way forward at the {\em pi\`u animato} with a change of tempo and the introduction of new material.  Wagner's resolution is inconclusive, in line with the shifting sensibilities of the protagonists. 
A glance at the Levi graph in Figure \ref{eulerian_tonnetz} is revealing. One sees that Chopin's resolution takes the initial $G_{\sharp}m$ along the path of a 
$2p$-decacycle to its polar $B_{\flat} M$. From $G_{\sharp}m$ we follow the ``floppy bow tie'' trajectory 
to $BM$, then to $E_{\flat} m$, then across to $E_{\flat}M$, up to $Gm$, and finally to $B_{\flat}M$ at the opposite side of the decagon.
It seems that it is arrival at the polar that gives this resolution its firmness. 
The situation with Wagner's ambivalent resolution is  less decisive: 
he moves the chord to an adjacent one, one step clockwise along the unique 3$p$-hexacycle to which both chords belong  from $G_{\sharp}m$ to $EM$. The simplicity of this transformation suggests why its character is so unsettled: nothing has  changed. One is moving sideways, from indecisiveness to indecisiveness. 
Indeed, by studying the positions of the reductive triads on the Levi graph and how the one chord can be reached from the other along a trajectory of a cycle of the tonnetz, one is offered insight into the musicality of such transitions.

There is a substantial literature addressing  how tetrachords can be accommodated into tonnetz-like schemes (Childs 1998, Douthett and Steinbach 1998, Cohn 2012, Tymoczko 2011, Nu\~no 2021). Such schemes have  included dominant sevenths, half-diminished sevenths,  minor sevenths, diminished sevenths and French sixths. 
We present our own approach  in Section VI in the form of a new tonnetz for Tristan-genus chords.  For the moment we take advantage of the fact that tetradic relations can be given an illuminating analysis in terms of the consonant triads underpinning them. See Cohn (2012) in Chapter 7, ``Dissonance'', at the section ``Reduction to a Triadic Subset'' (pp. 142-145); NB Cohn's remark at page 144, ``Wagner often treats minor triads and their $\varnothing^7$ supersets as interchangeable in his late music.'' See also Bailey (1985) at pages 122-125. 

\begin{figure} [htbp]
%\centering
\includegraphics[clip, scale=0.70]{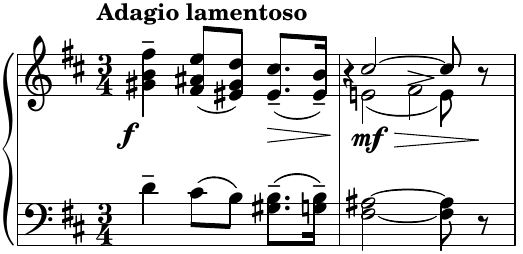}
\caption{Tchaikovsky, Symphony No.~6 in B Minor, Op.~74, fourth movement, opening bars. The initial  $Bm^6$  in the permutation $[D, G_\sharp, B, F_\sharp ]$ resolves to $F_\sharp ^{7}$. Transposed down by three semitones, the initial chord takes the form  
of a permuted $G_{\sharp}m^6$, namely  $[B, E_\sharp, G_\sharp, D_\sharp]$ or equally $[B, F, A_\flat, E_\flat]$, obtained from the Tristan chord by the permutation $[1,2,3,4] \to [2,1,4,3]$. The resolution is to $E_{\flat}^{7}$, in contrast with Chopin's resolution to $B_{\flat}^{7}$ and Wagner's to $E^{7}$.} 
\label{tchaikovsky}
\end{figure}

An interesting and rather distinctive example of the resolution of a Tristan-genus chord  can be found in the fourth movement (Adagio Lamentoso) of Tchaikovsky's Symphony No.~6 in B Minor, Op.~74 (``Path\'etique'').  Tchaikovsky's modified version of the Tristan chord is heard  in the strings at the outset of the Adagio Lamentoso. 
Recall that Wagner's chord is $[F, B, E_{\flat}, A_{\flat}]$, whereas Tchaikovsky, in line with the overall key of his symphony uses a $Bm^{6}$ chord, in the permutation $[D, G_{\sharp}, B, F_{\sharp}]$, as in Figure \ref{tchaikovsky}, resolving to ${F_{\sharp}}^{7}$. 

%
%%%%%%%%%%Tchaikovsky excerpt
%%%%%

If we transpose this passage down three semitones to facilitate comparison with Wagner and Chopin, the result is a modified Tristan chord in the form
$[B,  F, A_{\flat}, E_{\flat}]$, which resolves to ${E_{\flat}}^{7}$. 
Thus, Tchaikovsky swaps the two lower voices and the two upper voices of the Tristan chord and the resulting resolution of  $G_{\sharp}m^6$ is more like Chopin's than Wagner's, only leading to ${E_{\flat}}^{7}$ rather than ${B_{\flat}}^{7}$. 

Hence, as in the case of Chopin, there is a sense of definition in Tchaikovsky's resolution; but there is no joy in it -- this resolution has more the character of a grim facing up to reality with reluctant finality. We are speaking of the sense of emotion that Tchaikovsky creates so successfully in the opening measures of this movement. 
But how is it achieved? Any assignment of a particular musical structure as being evocative of a specific emotion has to be taken as speculative. The extent to which psychological notions can be meaningfully incorporated into the theory of music in a mathematically systematic way remains to be seen. Clearly, even the idea of perceived  distance between chords is not without difficulty. Tymoczko (2011) offers a critique of tonnetz-based distance measures in his Appendix C. Nonetheless, when such identifications can be made it is worth pointing them out. In the present context, it is by means of a $2p$-hexacycle. 
Specifically, beginning at $G_{\sharp}m$ in the guise of a permuted $G_{\sharp}m^{6}$ (we stick with Wagner's pitch classes), we find there is a ``straight bow tie'' $2p$-hexacycle that takes one to $A_{\flat}M$ then up to $Cm$ then finally to $E_{\flat}M$. 
As in Chopin's resolution, following the hexacycle right through to the polar position contributes to the sense of completion.  
However, Tchaikovsky's resolution to $E_{\flat}M$, since it is only across a $2p$-hexacycle rather than a $2p$-decacycle, as in the case of Chopin,  thus has a more perfunctory character, that of a plagal cadence, and hence its desolation. 
%%%%%%%%%Parsifal excerpt
\begin{figure}[htbp]
%\centering
\includegraphics[clip, scale=0.70]{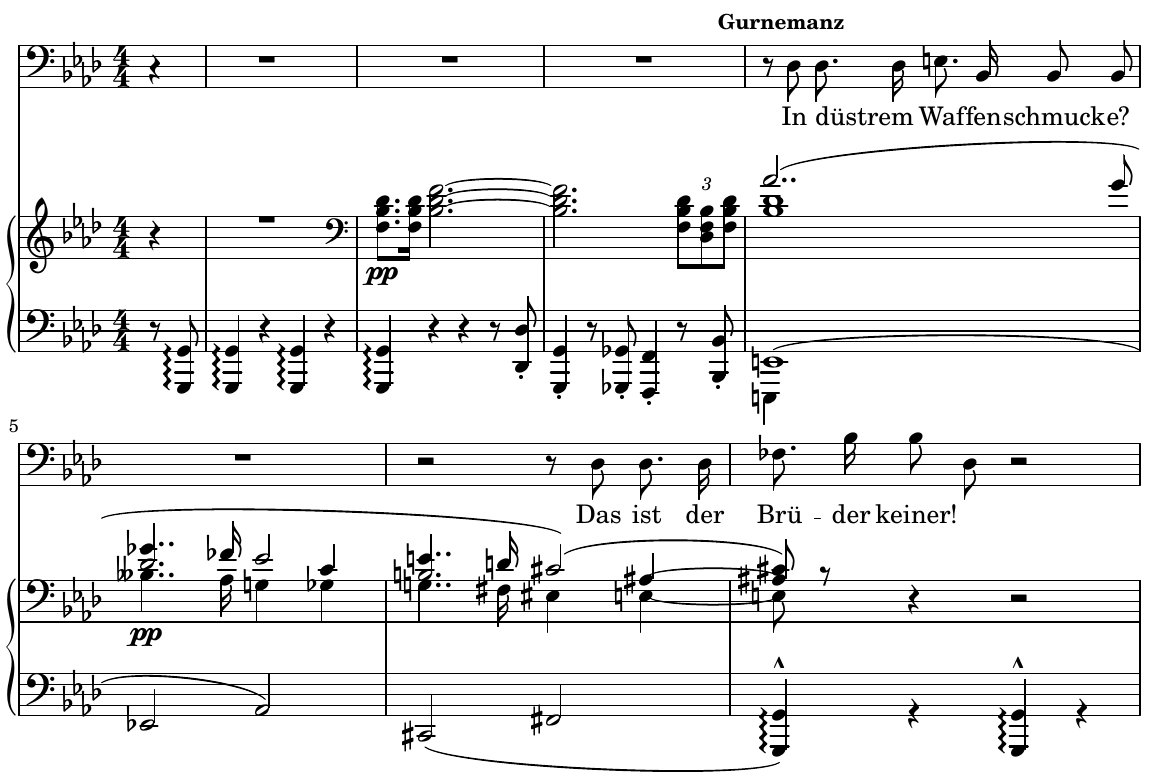}
\caption{Wagner, {\em Parsifal}, excerpt from Act 3. The Parsifal leitmotif, which at Parsifal's entry in Act 1 is a simple fanfare  leading from $B_{\flat}M$ to $E_{\flat}M$, has evolved at Parsifal's Act 3 entry into $B_{\flat}m^6$  leading to $D_{\flat}m^6$, which Wagner presents in the permutation $\{ E_{\natural}, B_{\flat}, D_{\flat}, A_{\flat}\}$, isomorphic to Tchaikovsky's $Bm^{6}$ in the form $\{D, G_{\sharp}, B, F_{\sharp}\}$ appearing at the outset of the Adagio Lamentoso. } 
\label{parsifal entrance}
\end{figure}
%%%%%

Did Wagner himself make use of Tchaikovsky's version of the Tristan chord? The reply is ``yes.'' Tchaikovsky's chord can be found in {\em Parsifal}, early in Act 3, at the passage where the horns announce  the return of Parsifal to the now languishing Kingdom of the Holy Grail, disguised in black armour and carrying the Sacred Spear. 

The Parsifal motive, which in the first act took the form of a simple fanfare involving a triplet figure leading from $B_{\flat}M$ to  $E_{\flat}M$ has morphed  into a fanfare on $B_{\flat}m^6$  leading to a devastating $D_{\flat}m^6$ (Figure \ref{parsifal entrance}). 
It is this $D_{\flat}m^6$, which Wagner presents in the form $[ E_{\natural}, B_{\flat}, D_{\flat}, A_{\flat}]$, which is isomorphic to the $Bm^{6}$ that Tchaikovsky presents in the form $[D, G_{\sharp}, B, F_{\sharp}]$ in the Adagio Lamentoso. 
Wagner's version of the chord is pitched a tone higher than that of Tchaikovsky, but lowered by an octave with the 
$E_{\natural}$ doubled in the bass to add darkness and weight.  
It does not resolve to a dominant seventh, but rather, and only tentatively, to  a $G_{\flat}m^6$ in the form $[ E_{\flat}, B_{\flat\flat}, D_{\flat}, G_{\flat} ]$, which is a Tristan chord. 

Despite the surrounding sense of gloom and oppression, Wagner's resolution is ultimately conventional, in keeping with Parsifal's character, that is, up a fourth from $D_{\flat}m^6$, in line with the basic Parsifal motive -- and hence consoling. 
The Tristan chord then resolves further, with more certainty, to another Tristan chord a tone lower given by an $Em^6$ in the form $[ C_{\sharp}, G, B, E]$. 
Finally, the resolution proceeds to a diminished seventh $[G, E, A_{\sharp}, C_{\sharp}]$ reflecting Gurnemanz's bewilderment at the arrival of the strange knight (``Das ist der Br\"uder keiner!''). 
Wagner's resolution of the ``Tchaikovsky chord'' can be best understood in the context of the phrase leading from the horn fanfare in 
$B_{\flat}m^6$ to the Tchaikovsky chord at $D_{\flat}m^6$ on to the first Tristan resolution at 
$G_{\flat}m^6$ and then to the second Tristan resolution at $Em^6$. The starting point  $B_{\flat}m^6$ and the end point $Em^6$ are polar on the Hamiltonian cycle of the tonnetz graph shown in Figure \ref{eulerian_tonnetz}. 

The Tchaikovsky chord is one quarter around the tonnetz, half way between the starting point and the end point. The first Tristan resolution is one black vertex further along on the cycle towards the end point, and can be viewed as a staging point before the final stage is reached at the second Tristan resolution. 
This brief passage is highly symmetrical in its construction, spanning the full breadth of the tonnetz, even though with its gritty dissonances and abrupt conclusion, and only lasting a few seconds, it has the character of a throwaway line. We observe that all four of the chords are linked by a 4$p$-decacycle, 
\begin{eqnarray}
  \langle B_{\flat}m,  G_{\flat}M, G_{\flat}m, AM, D_{\flat}m, EM, Em, GM, Gm, B_{\flat}M, B_{\flat}m \rangle.
\end{eqnarray}
The diminished seventh $[G, E, A_{\sharp}, C_{\sharp}]$ that finishes off this sequence of four minor sixths representing Gurnemanz's take on the unexpected guest is merely one semitone away from a $Gm^6$, and hence for Wagner's dramatic purposes serves well as a point of exit from the decacycle. This view is not inconsistent with the idea that the Hamiltonian cycle circumnavigating the tonnetz is also playing a role -- for not only is $Em$  opposite to $B_{\flat}m$ in the Hamiltonian cycle, one also finds that $Gm$ is opposite to $D_{\flat}m$. Hence the chords $B_{\flat}m$, $D_{\flat}m$, $Em$, $Gm$ are equally spaced around the perimeter of the tonnetz and the names of these chords spell out the notes of the diminished seventh on which the sequence ends. 

Tchaikovsky's chord can also be found at the strident opening of Chopin's Scherzo in B Minor, Op.~20, in the form $[G, C_{\sharp}, E, B]$, resolving to $[A_{\sharp},  C_{\sharp}, E, F_{\sharp}]$, both with some doubling.  Thus, we have $Em^6$ moving to ${F_{\sharp}}^{7}$, the first chord being in Tchaikovsky's permutation and the second being, unusually, in the first inversion. Putting this into Wagner's Tristan pitch levels, we have $G_{\sharp}m^6$ moving to ${A_{\sharp}}^{7}$, which we recognize as enharmonically equivalent to the resolution in the Ballade in G Minor at the {\em pi\`u animato}; though in the Scherzo, since the dominant seventh appears in the first inversion, it has a different character, demanding immediate further resolution to the tonic $D_{\sharp}m$, which in the key of the Scherzo is $Bm$. 
It is interesting to note that at Act 1, Scene 5, mm.~6-7 of \emph{Tristan und Isolde} Wagner resolves the Tristan chord in its original pitches to ${B_{\flat}}^{\!7}$.  A direct analogue of the descending figure in the right hand in the two bars preceding the \emph{pi\`u animato} of Chopin's Ballade can be found at the climax of the Tristan Prelude, in the violins, following the fortissimo Tristan chord two bars before the \emph{Il tempo poco a poco ritenuto}.

%%%%%%%%%%%%%%%%%

\section{Tonnetze and Tessellations}
\label{sec:Tone Networks and Tessellations}

\begin{figure} [htbp] %Archimedean tonnetz
\centering
\includegraphics[clip,scale=0.50]{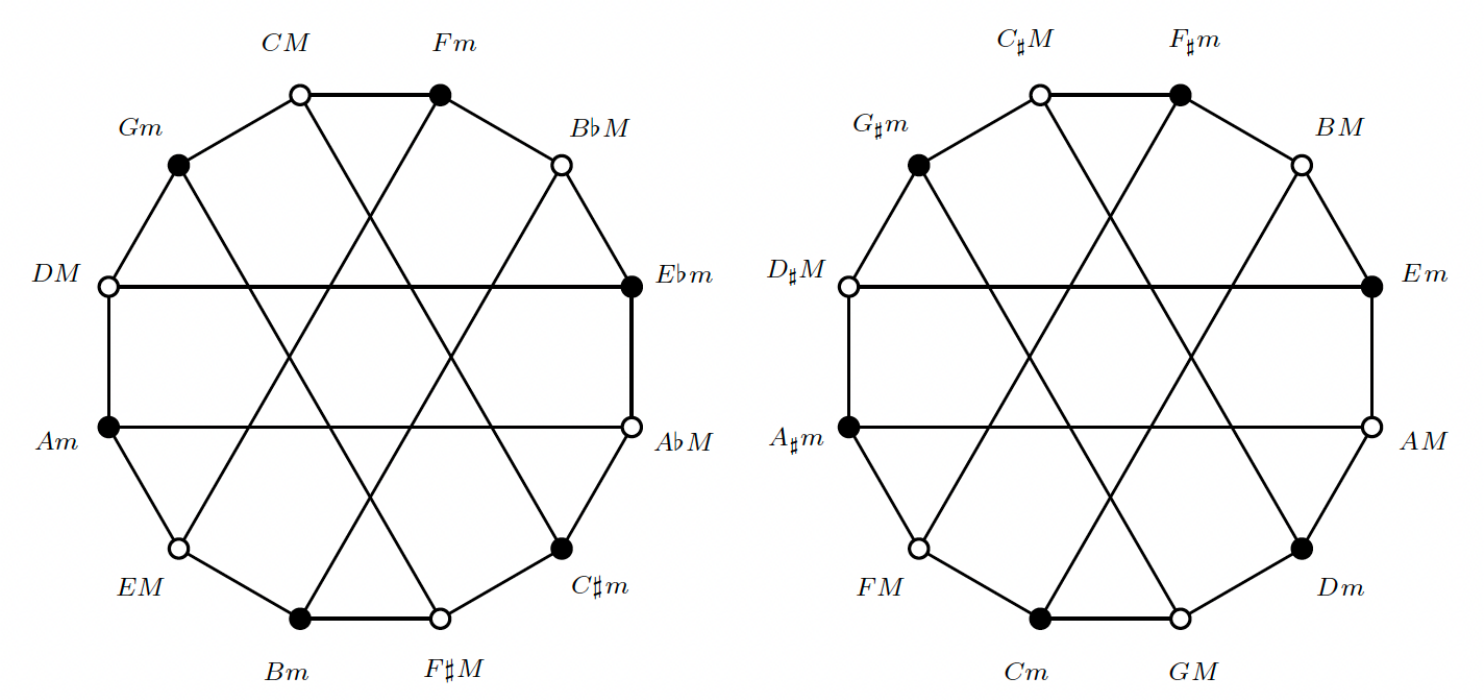}
\caption{The Archimedean tonnetze. Each major triad can be changed into a minor triad by shifting two of the tones, and each minor triad can be changed into a major triad by shifting two of the tones. The resulting bipartite graph has two connected components.}
\label{Archimedean tonnetz}
\end{figure}

\noindent As a first step towards construction of a tonnetz for tetrachords, we give another example of a tone network that can be constructed by purely combinatorial arguments and yet admits a geometric representation that is not unpleasing.
Given a major triad one can ask for all the minor triads that can be constructed from it by holding one tone fixed and altering the other two. Thus, given $CM$ we can construct three minor triads differing from it at two of the tones -- namely, $Fm$, $Gm$ and $C_{\sharp} m$. 
Similarly, given $Cm$, we can construct  three major triads differing from it at two of the tones -- namely, $FM$, $GM$ and $BM$ \footnote{The single-common-tone preserving relations between triads that we have described are known as the ``obverse'' relations of the usual $\rm\bf L$, $\rm\bf  P$ and $\rm\bf  R$ operations. For example, the map taking $CM \to Fm$ is $\rm\bf  L'$, the obverse of $\rm\bf  L$; $CM \to C_{\sharp} m$ is $\rm\bf  P'$, the obverse of $\rm\bf  P$; $CM \to Gm$ is $\rm\bf   R'$, the obverse of $\rm\bf  R$. These terms were coined by Morris (1998). Note that $\rm\bf L'$  = $\rm\bf RLP$, $\rm\bf P'$ = $\rm\bf RPL$, and $\rm\bf R'$ = $\rm\bf LRP $ (Cohn 1997). Our characterization of the tonnetz generated by these three maps and their inverses as a bipartite cubic graph of girth four in two components is new, as far as we are aware, as is its representation as a tessellation.}. The associated bipartite graph is not connected, but rather is given by the union of two distinct connected components. Thus we obtain two sets of twelve triads, each containing six major triads and six minor triads, as shown in Figure \ref{Archimedean tonnetz}. 

%%%%%%%%%%%%%%%%%%%%%%%

\begin{figure} [htbp]
\includegraphics[scale=1.15]{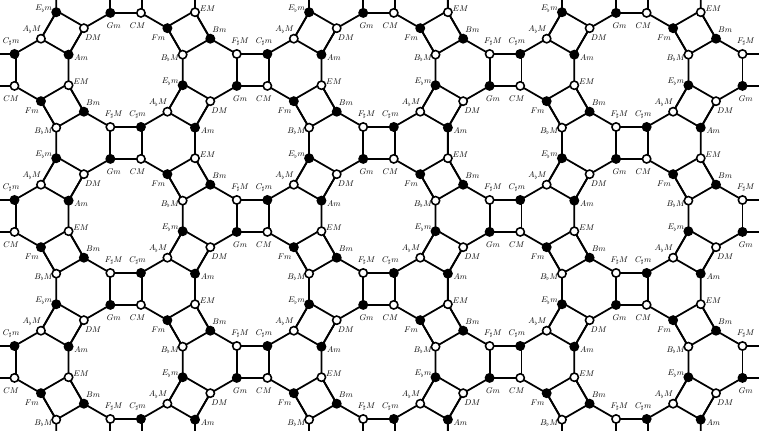}
\caption{Tessellation of the plane with dodecagons, hexagons and squares generated by the first of the Archimedean tonnetze. A similar tessellation can be constructed for the second one.}
\label{dodecagontessellation}
\end{figure}

As with the Eulerian tonnetz, one can read off the resulting cycles directly from the bipartite graph (see Appendix). There are hexacycles, octacycles, decacycles and dodecacycles, but  also tetracycles like $\langle Gm$, $F_{\sharp}M$, $C_{\sharp}m$, $CM$, $Gm\rangle$. 
The significance of the tetracycles is that these graphs are not the graphs of configurations. For a bipartite graph to be a Levi graph  it must be of girth at least six.
Nonetheless, if we lay out copies of the dodecagons in a lattice and connect adjacent vertices across the lattice (rather than across the dodecagon), what emerges is a striking tessellation of the plane, composed of dodecagons, hexagons, and squares, as in Figure \ref{dodecagontessellation}.
%%%%%
%%%%%
This tessellation is one of a number known to Kepler (1619). A dual tessellation called the ``Laves tiling'' associated with $\{4, 6, 12\}$ can be constructed (Gr\"unbaum and Shephard 2016, pages  95-98) with the property that the vertices of the original tiling correspond to the tiles of the dual tiling and the tiles of the original tiling correspond to vertices of the dual tiling. The edges of the original tiling and its dual can then be identified. The Laves tilings associated with the Archimedean tonnetze can be seen in Rietsch (2024), pages  336-338. Each dodecagon of the original tonnetz corresponds to a vertex of the corresponding Laves tessellation met by twelve edges. Thus, our geometric/combinatorial method gives a construction of the duals of Rietsch's $G_2$ tonnetze, which were obtained with Lie theory.  
Each dodecagon is surrounded by six hexagons and six squares in a remarkable pattern known as an Archimedean or semi-regular tessellation of type $\{4, 6, 12\}$. The internal angle of two adjacent sides of a regular $n$-gon is 
$(n-2)\pi/n$. It follows that if we wish to fit three polygons all with sides of the same length around the same vertex, then the numbers of sides $n_1, n_2, n_3$ of the polygons must be so that $(n_1-2)\pi/n_1 + (n_2-2)\pi/n_2 +(n_3-2)\pi/n_3 = 2\pi$. This condition is satisfied by $n_1 = 4$, $n_2 = 6$, $n_3 = 12$. More generally, there are 21 different solutions with various numbers of polygons and various numbers of sides to fit a collection of regular polygons around a vertex. Of these, eleven can be extended across the plane, giving the eleven  Archimedean tilings. 
See, for example, Gr\"unbaum and Shephard (1977, 2016), Conway, Burgiel and Goodman-Strauss (2008) or Wilson (2016). 

Since the triads at adjacent vertices share a tone in common, we can associate each edge with a tone. The three edges meeting a vertex give the tones of which the triad at that vertex is composed. The Hamiltonian cycles in the Archimedean tonnetze ascend in the circle of fifths going counterclockwise around the rim. As a result, two of the three triads a triad is adjacent to are a fourth or fifth apart.
The Archimedean tonnetze can be used as a basis for composition and they are rather useful in this respect. 
The fact that two tones of a triad change at each transition ensures that the chord sequences are less bland than the hexacycles of the Eulerian tonnetz, and hence are more interesting to listen to even in the absence of the rich tonal environment from which historical examples of the use of the Eulerian tonnetz are usually extracted. As Piston (1985) says, ``The purpose of this procedure [of voice leading in harmonic progressions] is to insure the smoothest possible connection of two chords, so that one seems to flow into the next. Continued practice of the process will, however, result in rather dull music.''

\begin{figure}  [htbp]
 \centering
\includegraphics[scale=0.50]{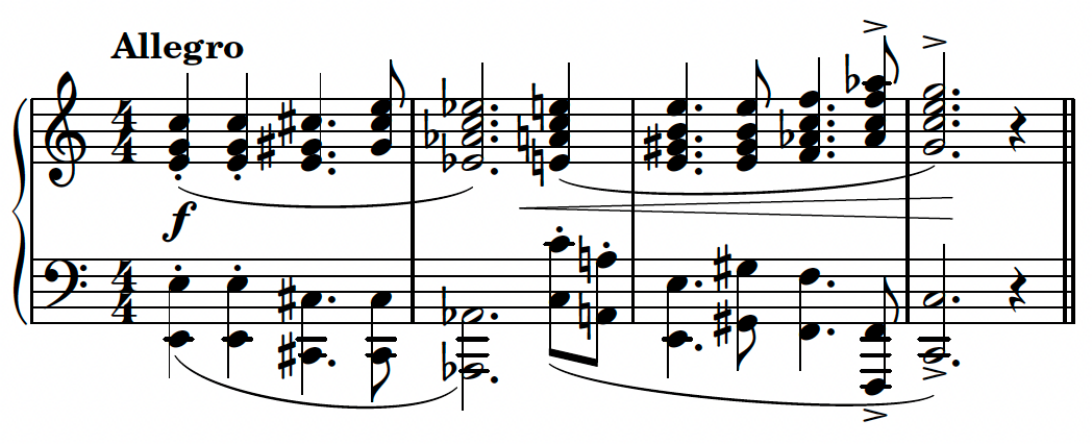}
\caption{A simple ``Fanfare for Kepler'' based on a hexacycle of the first Archimedean tonnetz. The cycle $\langle CM, C_{\sharp}m, A_{\flat}M, Am, EM, Fm, CM \rangle$ used here is a closed sequence of six major and minor triads with the property that each triad differs from its successor at exactly two tones.}
\label{fanfare}
\end{figure}

In Figure \ref{fanfare} we illustrate this point with a short composition in the form of a fanfare on a cycle of the Archimedean tonnetz. 
The occasional common practice modalities embedded in this tonnetz, such as the plagal cadence from $Fm$  to $CM$, help to ground the piece with an aura of tonality. The various cycles of the Archimedean tonnetz have contrasting voice-leading properties (cf.\,\,Tymoczko 2020, pages 116-128). A counterclockwise movement through the tetracycle $\langle CM, C_{\sharp}m, F_{\sharp}M, Gm, CM \rangle$ in Figure \ref{Archimedean tonnetz} induces a transition from $CM^{(0)}$ to $CM^{(1)}$, whereas movement through a hexacycle returns each voice to its starting pitch and a clockwise trajectory through the 0$p$-dodecacycle leads every voice up an octave. Thus, the figure reflects aspects of the topology of voice-leading spaces. The tessellation in Figure \ref{dodecagontessellation} has an advantage over the graph of Figure \ref{Archimedean tonnetz} since the various cycles are visually more distinct.
%

%%%%%%%%%%%%%%%%%%%
\section{Tone Networks for Tetrachords}
\label{sec:Tone Networks for Tetrachords}

\begin{quote}`` .\,.\,.\,indeed, the {\em Tristan} prelude (and the opera as a whole) could be said to be ``about" the various ways of resolving a Tristan chord to the dominant seventh chord\,.\,.\,.''

\hfill ---D.~Tymoczko,~{\em A Geometry of Music}
\end{quote}
%
%\vspace{0.25cm}

\noindent The investigations that have been hitherto carried out on the construction of tone networks for tetrachords offer many insights into the relation of the families of tetrachords to one another but at the same time one is forced to realize the daunting nature of the task. In {\it Tristan und Isolde}, the enharmonic minor sixth chord (that is, the Tristan chord in one of its permutations)  resolves in no fewer than eight distinct ways into dominant seventh chords, according to Tymoczko (2011). If we transpose the initial chord to $G_{\sharp}m^6$, the eight resolutions starting from that chord are as follows:
(i)$ \hspace{0.15cm} G_\sharp m^6 \to E^7$, 
(ii)$ \hspace{0.15cm} G_\sharp m^6 \to F^7$,
(iii)$ \hspace{0.15cm} G_\sharp m^6 \to G^7$, 
(iv)$ \hspace{0.15cm} G_\sharp m^6 \to  {A_\flat}^{\! 7}$,
(v) $\hspace{0.15cm} G_\sharp m^6 \to {B_\flat}^{\! 7}$,
(vi)$ \hspace{0.15cm} G_\sharp m^6 \to B^7$, 
(vii)$ \hspace{0.15cm}  G_\sharp m^6 \to {D_\flat}^{\! 7}$, and
(viii)$ \hspace{0.15cm} G_\sharp m^6 \to D^7$.
This list includes Chopin's resolution to ${B_\flat}^{\! 7}$ in the $G$-minor Ballade, but not Tchaikovsky's resolution (after transposition) to ${E_\flat}^{\!7}$ in the Adagio Lamentoso, which nonetheless works quite well as a resolution. Perhaps it too ought to be added to the list on an honorary basis? It is tempting to believe  that almost any resolution to a dominant seventh is acceptable. Are the ones listed above merely the ones that Wagner happened to get around to using? Is there any logic to Wagner's choice of admissible resolutions? 

Our approach to the construction of a tonnetz for minor sixths and dominant sevenths will be to adopt a variant of the strategy that we used to create the Archimedean tonnetze. Thus, we consider all transitions from minor sixth chords to dominant sevenths with the property that exactly two tones are altered. 
The goal is then to derive a Levi graph for the resulting relations. A problem that immediately arises is that there are too many such transitions, some of them even admitting tetracycles, such as 
$\langle {G_\sharp m}^{6}, B^7, {Dm}^{6}, F^7, {G_\sharp m}^{6} \rangle$, 
hence blocking the possibility of a Levi graph. Thus, we need some principle that can be applied uniformly across the class of chords being considered to limit the range of transitions in a way that is sufficient to ensure that we get a Levi graph. Our criterion for success will be that the resulting Levi graph is associated with an interesting configuration.

\begin{figure} [htbp]
\includegraphics[scale=0.70]{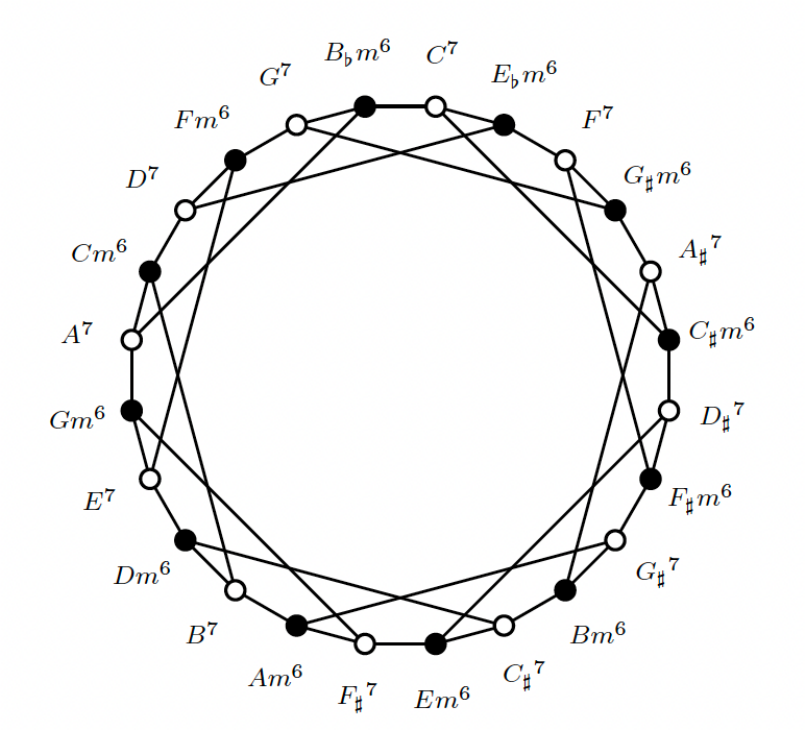}
\caption{A Tristan-genus tonnetz. The twelve minor sixth chords and twelve dominant seventh chords can be arranged in a Levi graph with the incidence geometry of a  configuration $\{12_3\}$. Each minor sixth chord is connected to three dominant seventh chords and each dominant seventh is connected to three minor sixths. At each transformation two notes are retained in the new chord, one being the ``dissonant'' note of the original tetrachord. Edges that cross the tonnetz preserve the dissonant note of the tetrachord as a dissonance in the new tetrachord.}
\label{fig:tristan_tonnetz}
\end{figure}

Classical music theorists have observed that in both the minor sixth and in the dominant seventh there is a single ``dissonant'' tone that combines with a consonant triad. For example, the dissonant note in $Cm^6$  is  $A$ and in $C^7$ it is $B_{\flat}$. Thus, $Cm^6$ is the $C$-minor triad ``with an added sixth,'' and $C^7$ is the $C$-major triad ``with an added minor seventh,'' the ``added'' notes being the dissonant ones.
Let us therefore impose the condition that the dissonant note is unchanged. That implies that two members of the underlying consonant triad shift, and we are back to a situation similar to what we considered in Section \ref{sec:Tone Networks and Tessellations} in connection with the Archimedean tonnetze. 

The result is constrained by the requirement that a minor sixth must transform into a dominant seventh and vice versa. 
Although the dissonant tone in the tetrad being mapped is preserved, it need not necessarily map to a dissonant tone of the target.  And likewise, a consonant note of the initial triad might map to a dissonant tone of the target. 

It is easy to see that there are exactly three two-note transitions that send $Cm^6$ to a dominant seventh chord in such a way that the note $A$ is retained, namely $Cm^6 \to A^7$, $C m^6 \to B^7$ and $Cm^6 \to D^7$.
Likewise, there are exactly three ways in which $C^7$ can be sent to a minor sixth by a two-note transition such that the note $B_\flat$ is retained, namely $C^7 \to  B_{\flat}m^6$, $C^7 \to E_{\flat}m^6$ and $C^7 \to C_{\sharp}m^6$.
Similar transformations exist by transposition for the other minor sixths. Putting these together, we obtain a connected bipartite graph of degree three and girth six for the minor sixths and dominant sevenths, shown in Figure \ref{fig:tristan_tonnetz}. We thus obtain a Levi graph and a configuration
$\{12_3\}$. 

By counting cycles, one easily determines that this graph is different from the Levi graph of the Eulerian tonnetz.  Now, the Eulerian tonnetz has a total of sixteen hexacycles, including the four ``hexatonic'' 3$p$-hexacycles and the twelve ``pitch-retaining'' 2$p$-hexacycles. But the Tristan-genus tonnetz of Figure \ref{fig:tristan_tonnetz} admits twenty-four hexacycles -- specifically, we have twelve 1$p$-hexacycles (hexabeanies) of the form 
$\langle C^7, \,E_\flat m^6, \,F^7, \,G_\sharp m^6, \,A_\sharp^7, \,C_\sharp m^6, \,C^7 \rangle$
and twelve 2$p$-hexacycles (bow ties) of the form
 $\langle C^7, \,E_\flat m^6, \,F^7, \,F_\sharp m^6, \,D_{\sharp}^{7}, \,C_\sharp m^6, \,C^7 \rangle$.
 
 For convenience, we have retained here the ``$p$'' notation for edges that cross the graph, where in the present context one can take it to mean ``preserve'' (as in ``preserving the dissonant tone'') rather than ``parallel''. 

\begin{figure}[htbp]
\includegraphics[scale=0.55]{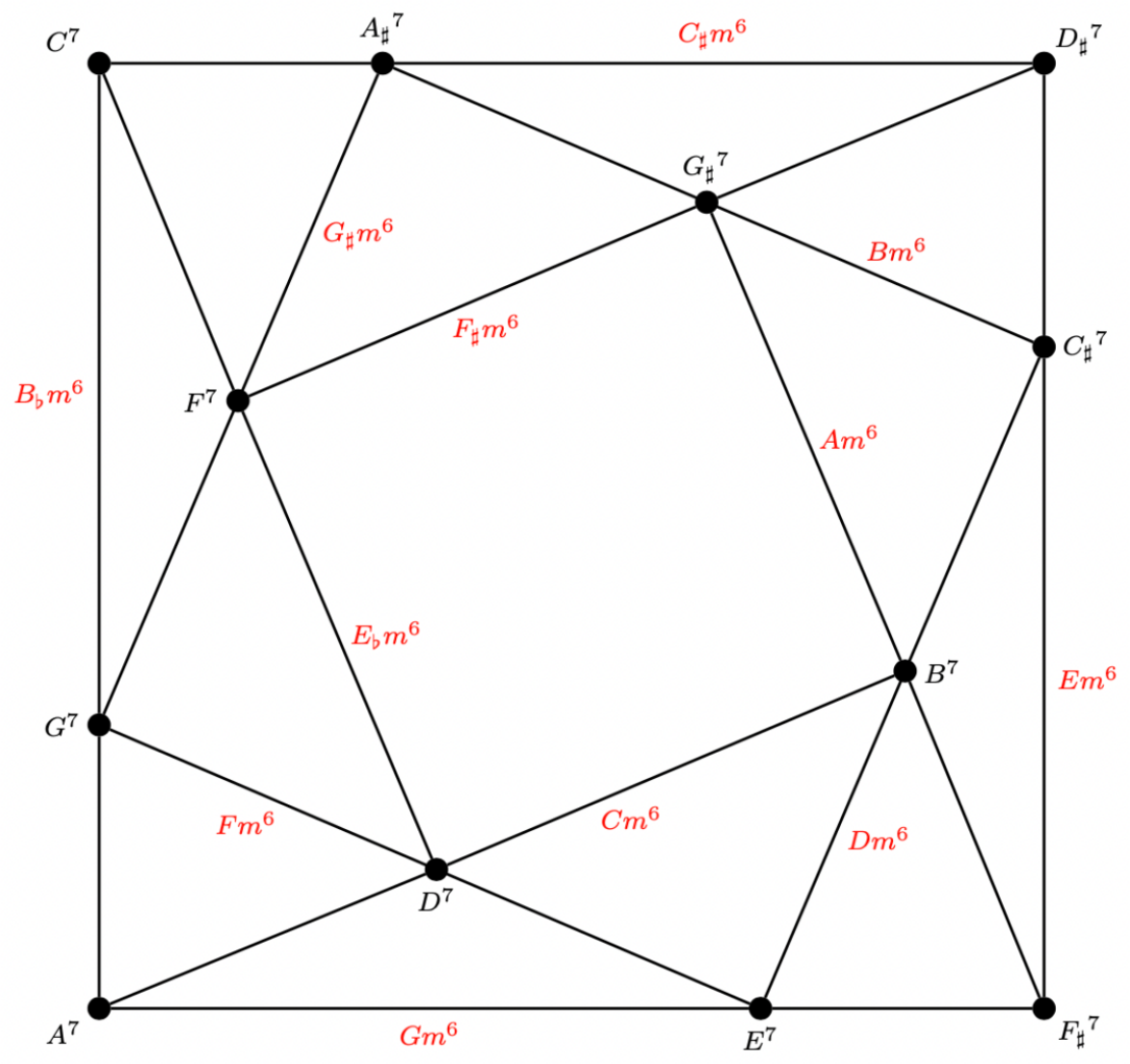}
\vspace{0cm}
\caption{The geometric configuration of the Tristan-genus tonnetz is the D228 of Daublebsky von Sterneck, here depicted. Its Levi graph appears in Figure \ref{fig:tristan_tonnetz}. This self-dual configuration is one of the four ``homogeneous'' cases arising among the 229 configurations of type $\{12_3\}$, alongside D222, D226, and D88. The group of the configuration is generated by quarter-turn rotations of the overall figure and cyclic permutations of the three mutually inscribed squares among each other. The squares correspond to the three $4p$-octacycles of the Levi graph.}
\label{tristan_tonnetz_configuration}
\end{figure}

Since our new Levi graph has a different number of hexacycles than that of the Eulerian tonnetz, this means that the two are not isomorphic. 
On the other hand, since it is known that all the combinatorial  configurations $\{12_3\}$ are realizable in the Euclidean plane, one can ask which of the 229 configurations of this type it turns out to be, if it is not the D222. 
We have been able to identify it as the D228 of Daublebsky von Sterneck (1895). This rather attractive and highly symmetric configuration is plotted in Figure \ref{tristan_tonnetz_configuration}.  

Since they are isomorphic, the Levi graph shown in Figure \ref{fig:tristan_tonnetz} and the configuration shown in Figure \ref{tristan_tonnetz_configuration} can each be referred to as a ``Tristan-genus tonnetz.'' 
Hence, we arrive at the following conclusion.

 %%%%%%%%% 
\begin{Proposition}
A tonnetz for the Tristan genus of dominant seventh chords and minor sixth chords can be constructed by limiting the maps between these chord sets to  two-tone transitions that preserve the dissonant tone of the chord being mapped. The resulting incidence geometry, which is distinct from that of the Eulerian tonnetz, is the self-dual configuration $\{12_3\}$  in $\mathbb {R}^{2}$ known as the {\em D228}  of Daublebsky von Sterneck. The $12$ dominant sevenths are represented by points and the $12$ minor sixths are represented by lines. The $36$ incidence relations between the $12$ points and $12$ lines determine the edges of the associated Levi graph, where the white vertices are dominant sevenths and the black vertices are minor sixths. 
\label{Tristan proposition}
\end{Proposition}
 
The symmetry group of this Tristan-genus configuration is generated by quarter-turn physical rotations of the overall figure and cyclic permutations of the three mutually inscribed squares among each other. The three squares, which are on an equal footing, correspond to the three $4p$-octacycles of the Levi graph. 
\begin{figure}[htbp]
\includegraphics[scale=0.60]{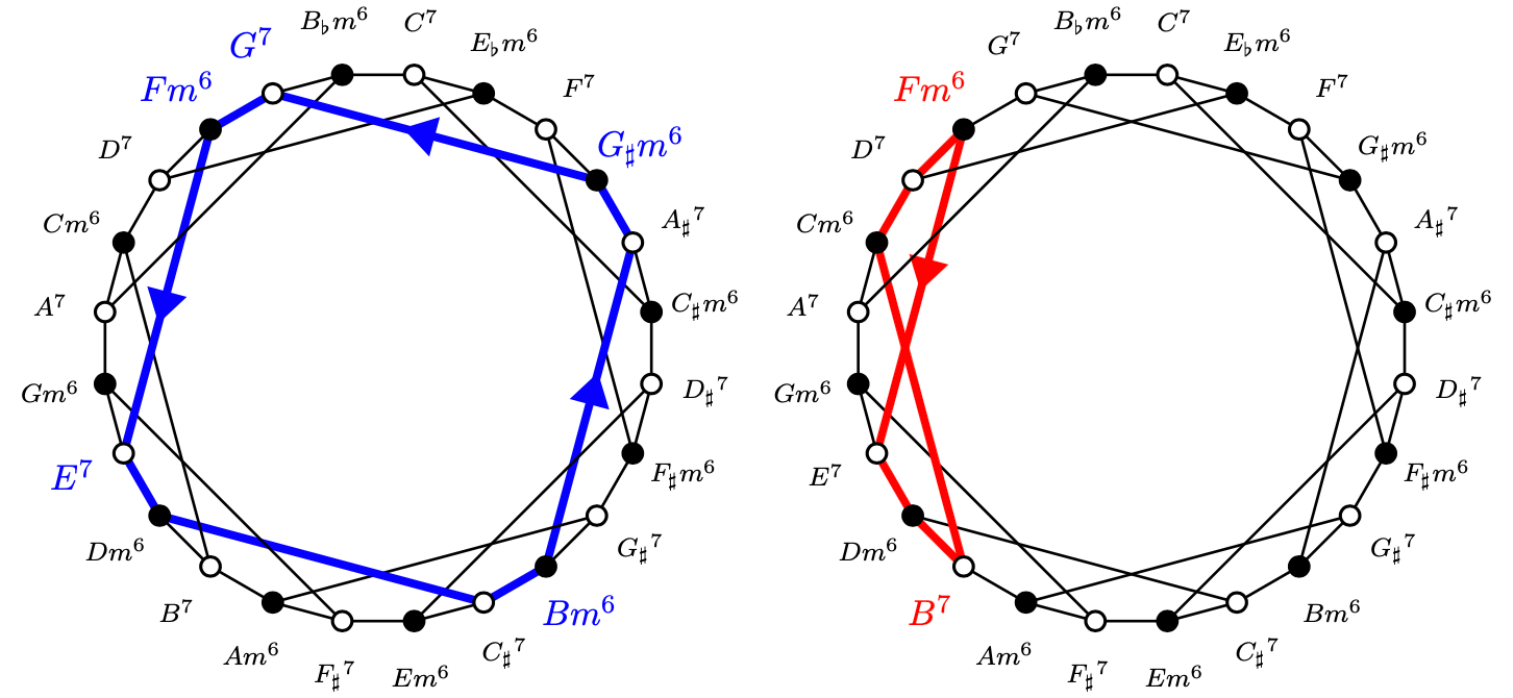}
\vspace{0cm}
\caption{Tonnetz analysis of the opening of the Tristan prelude. The twelve Tristan chords (each represented by a pitch-class equivalent minor sixth chord) lie in sets of four on the three $4p$-octacycles of the Tristan-genus tonnetz. The initial Tristan chord $G_\sharp m^6$ resolves to $E^7$ which lies opposite to it on the same octacycle, in the subpolar position. The second Tristan chord $Bm^6$ then resolves to its subpolar $G^7$, again on the same octacycle. The third Tristan chord  $Fm^6$ is also on the same octacycle, but resolves to $B^7$, which is polar to $Fm^6$ on the unique $2p$-hexacycle containing those two chords, and also acts as the second subpolar to the original $G_\sharp m^6$.}
\label{tristan_prelude_tetrachord_analysis}
\end{figure}
\vspace{0.20cm}

{\bf The Tristan Prelude}. The musical significance of the three mutually inscribed quadrilaterals of the configuration is that movements within these structures give rise to some of the characteristic sonorities  of late Wagner operas. 
We shall look at several examples of the use of these Wagnerian octacycles. 
First, we return to the Prelude to Act I of {\it Tristan und Isolde}. Three different versions of the opening motive 
appear in succession. At the first occurrence of the motive, shown in Figure \ref{fig:Tristan}, the Tristan chord $[F, B, E_{\flat}, A_{\flat}]$, which is pitch-class equivalent to $G_{\sharp} m^6$, resolves to $[E, G_{\sharp}, D, B]$, which is pitch-class equivalent to $E^7$.  
In its second occurrence, the motive is transposed up by a minor third, apart from the first tone, which is transposed up by a whole tone, so the interval between the first and the second tone is a major sixth rather than the minor sixth that one hears in the first occurrence of the motive. 
In the third occurrence, the initial rise is again by an interval of a major sixth, but the resolution of the minor sixth chord to a dominant seventh is different from that of the first two occurrences. In more detail, the resolutions are as follows:  
(a) $ G_\sharp m^6 \to E^7$, 
(b)  $Bm^6 \to G^7$,
(c)  $Fm^6 \to B^7$.
In an octacycle with the structure $\langle 1, 2, 3, 4, 5, 6, 7, 8, 1 \rangle$, where the odd numbers represent minor sixths and the even numbers represent dominant sevenths, we note that the polar of any chord is a chord of the same mode. Thus, the polar of chord 1 is chord 5, which is half way around the cycle. 
For some purposes we may be interested in the ``most opposite'' chords of the opposite mode. We can call these the ``subpolars''; there are two of them. Thus, the subpolars of chord 1 are chords 4 and 6. 
In the tables of cycle numbers shown in the Appendix we see that the Tristan tonnetz admits twelve 3$p$-octacycles and three 4$p$-octacycles. We observe that the polar of $G_\sharp m^6$ along the unique 4$p$-octacycle containing that chord is $Dm^6$. The subpolars of $G_\sharp m^6$ are $E^7$ and ${C_\sharp}^{\!7}$. 
Thus, Wagner's first resolution (a) runs counterclockwise along the octacycle beginning at $G_\sharp m^6$ and resolves at the first subpolar, $E^7$.  Figure   \ref{tristan_prelude_tetrachord_analysis} shows that resolution (b) begins on the same octacycle at $Bm^6$ and resolves to {\it its} first subpolar, which is $G^7$. 
Resolution (c) again starts on the same octacycle, this time at $Fm^6$, but resolves to $B^7$, which is {\it not} on the same octacycle -- but this $B^7$ nevertheless stands in a definite relation to $Fm^6$: in fact, $B^7$ is polar to $Fm^6$ along the 2$p$-hexacycle (a bow tie) that connects these two chords. 

In summary, as we see in the Levi graphs of Figure \ref{tristan_prelude_tetrachord_analysis}, all three of the minor sixth chords appearing at the outset of the prelude -- namely, $G_\sharp m^6$, $Bm^6$ and $Fm^6$ -- lie on the same octacycle, along with two of the resolvent dominant sevenths, $E^7$ and $G^7$. But to break off the cycle the final resolution starting at $Fm^6$ is deflected along one of the three 2$p$-hexacycles containing that chord, namely the one that passes through $E^7$ and  $Dm^6$, to reach its polar $B^7$, which we observe is the second subpolar of the original $G_\sharp m^6$. 
Hence, all six chords of this opening gesture of the prelude stand entwined in a rather subtle set of relations to one another on the Tristan-genus tonnetz  This passage at the beginning of the opera gives the impression of having sprung fully formed from the composer's mind. One is reminded of the remarkable analysis of William Blake's {\it Infant Sorrow} carried out by the linguist Roman Jakobson (1970), who after similar consideration of the intricate syntactic and morphological symmetries of Blake's short poem comes to the conclusion that it couldn't have been assembled piece by piece. 

\vspace{0.20cm}

{\bf Liebesnacht}. Another example of the use of octacycles in {\it Tristan und Isolde} can be found in the Liebesnacht duet of Act 2, Scene 2. This passage is marked by two warnings from Brang\"ane (``Habet acht!''). Following the first warning, Tristan pleads for death (``Lass mich sterben!''). Isolde responds similarly after the second warning. This section of the opera is accompanied by a family of six tetrachord progressions, each extending from minor sixth to minor sixth to dominant seventh to dominant seventh, visiting all twelve minor sixth chords along the way. The first of the six progressions is
\begin{eqnarray}
 A_\flat m^6 \to Dm^6 \to {D_\flat}^7  \to  C^7
 \label{Liebesnacht progression}
\end{eqnarray}
and the remaining five are transpositions of the first. The first three progressions, during which Tristan sings, begin respectively at $A_\flat m^6$, $B_\flat m^6$, and $Cm^6$, whereas the second three progressions, during which Isolde sings, begin at $Am^6$, $Bm^6$, and $C_\sharp m^6$.  In each progression the first three chords share a tritone. For example, in the progression \eqref{Liebesnacht progression} above the tritone common to $A_\flat m^6$, $Dm^6$, ${D_\flat}^7$ is the dyad
$\{F,B\}$. At the same time, the two minor sixth chords are themselves tritone-separated and hence polar on the Hamiltonian of the Tristan-genus tonnetz. They are also polar on the unique 4$p$-octacycle containing them. For instance, $A_{\flat}m^6$ and $Dm^6$ are polar on the perimeter Hamiltonian and on the 4$p$-octacycle 
$\langle A_{\flat}m^6, \,G^7,  \,Fm^6, \,E^7, \,Dm^6, \,{C_{\sharp}}^7, \,Bm^6, \,{A_{\sharp}}^7, \,A_{\flat}m^6 \rangle$.
In each of the six progressions, the first minor sixth chord and the final dominant seventh chord are simultaneously polar on the $2p$-hexacycle containing them and on the two $1p$-hexacycles containing them. Although it is tempting to dismiss the downshift in the seventh chords from ${D_\flat}^7$ to $C^7$ as a flourish, analogous to the big orchestral downshift from $AM$ to $A_{\flat}M$ towards the end of the Prelude to Act 3 of {\it Siegfried}, the way  the $C^7$ is cemented into the geometry here suggests otherwise.

\vspace{0.20cm}

{\bf Tymoczko's list}. Continuing with our analysis of {\it Tristan und Isolde}, we return to the list of the eight ``prominent Tristan-chord resolutions from the opera'' (Tymoczko 2011, Figure 8.6.10, page 302). Can we make sense of this list in the light of what we have learned? First, we observe that the names of these chords spell out the tones of the octatonic scale $O_{23}$. In particular,  the list splits into two sets, differing by a semitone: namely, $E^7, G^7,  {B_\flat}^{\! 7}, {D_\flat}^{\! 7}$, and  ${D}^{7}, F^7, {A_\flat}^{\! 7}, B^7$.  But these are the seventh chords that belong to precisely two of the three 4$p$-octacycles of the Tristan tonnetz. So now we understand which chords are missing from the list: they are the four dominant sevenths belonging to the excluded octacycle.  
Using the geometry of the configuration, one can prove a simple theorem, namely that once the Tristan chord has been selected, then the eight possible resolutions are determined. For there are three points of the configuration that lie on the line. Two of these lie at the vertices of one of the quadrilaterals, and the third will lie at a vertex of one of the other quadrilaterals. That fixes two quadrilaterals, and hence eight points of resolution. We have seen the eight points associated with the line 
marked $G_{\sharp}m^6$. To take another example at random, suppose we look at $B_{\flat}m^6$. This chord lies on the line through the three points $A^7$, $G^7$, and $C^7$. Thus the two quadrilaterals are those with the vertices $A^7$, $C^7$, $D_{\sharp}^7$, $F_{\sharp}^7$, and  $G^7$, $A_{\sharp}^7$, $C_{\sharp}^7$, $E^7$, and that gives us the eight resolutions. Now, you might say that surely one could simply transpose the eight resolutions of the original $G_{\sharp}m^6$ by the amount that $G_{\sharp}m^6$ itself differs from the new choice of the initial minor sixth. To be sure, this will work -- but the geometrical derivation is much quicker, and it is much less prone to error.

\begin{figure}[htbp]
\centering
\includegraphics[scale=0.50]{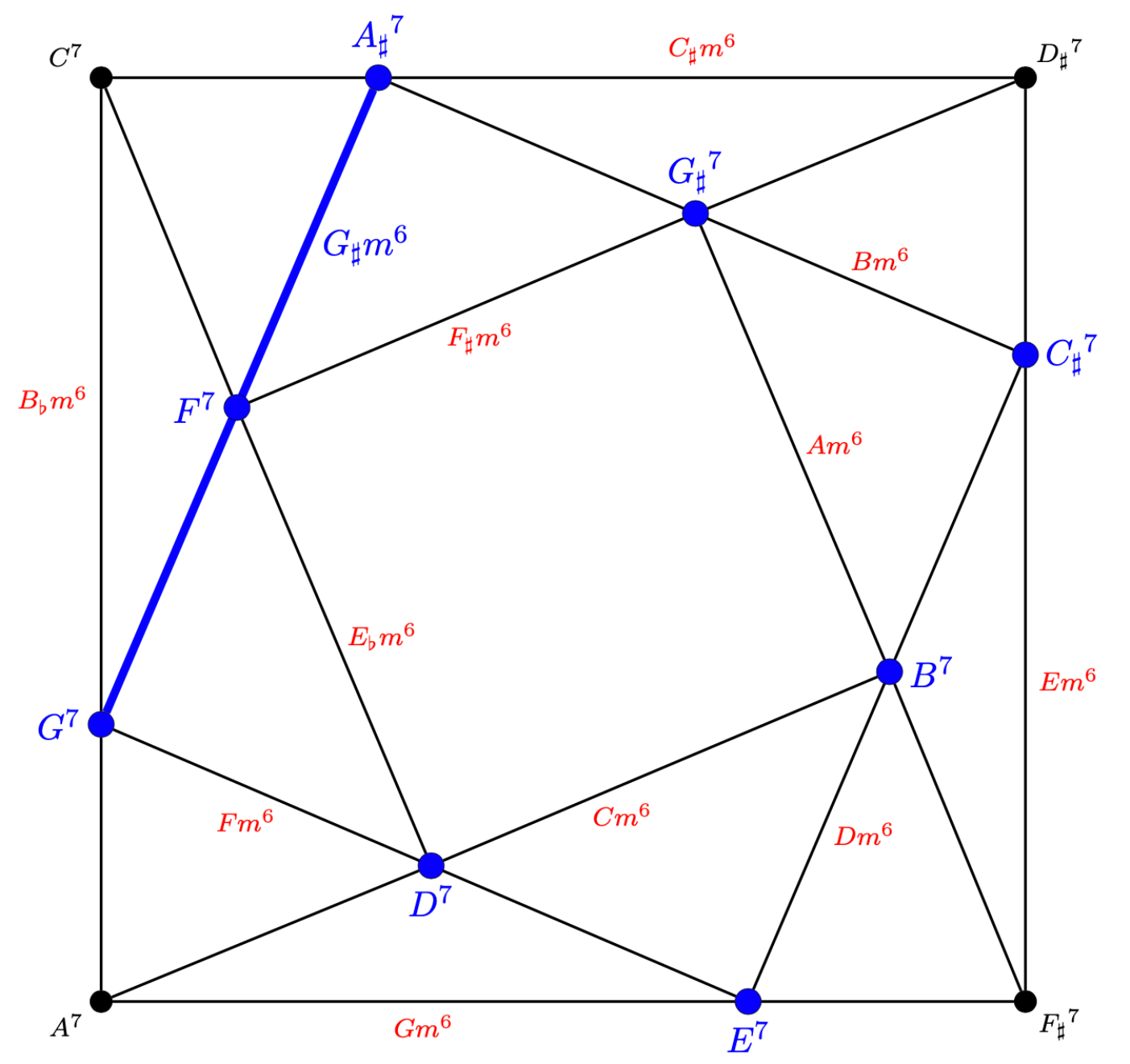}
\vspace{0cm}
\caption{Resolutions of the Tristan chord in the opera {\it Tristan und Isolde}. On the Tristan genus configuration, the eight resolutions occurring in the opera, when the initial chord is transposed to the line $G_{\sharp}m^6$, are marked by two sets of four points each, forming the vertices of a pair of quadrilaterals. The quadrilaterals represent octacycles, each consisting of an alternating sequence of four minor sixths and four dominant sevenths. The choice of the initial line (here, $G_{\sharp}m^6$) determines the two associated quadrilaterals, and hence the eight possible resolutions.} 
\label{Tymoczko_list}
\end{figure}
\vspace{0.20cm}

{\bf Tchaikovsky and Chopin}. How do the resolutions used by Tchaikovsky and Chopin that we had discussed in Section IV square up on the Tristan-genus tonnetz? 
In the case of Tchaikovsky, we observe that this resolution  from $G_\sharp m^6$ to $D_\sharp^7$ (in Wagner's pitches) at the beginning of the Adagio Lamentoso is polar along the unique 1$p$-hexacycle connecting these chords. In the case of the Eulerian tonnetz, the relation between the two underlying consonant triads is also via polarity across a hexacycle. 
In Chopin's resolution from $G_\sharp m^6$ to $B\flat^7$ at the transition to the {\em pi\`u animato} section of the Ballade in G Minor, we find that this $B\flat^7$ is immediately adjacent to $G_\sharp m^6$ on the Tristan-genus tonnetz. So, indeed, there is a kind of Tristan-genus tonnetz logic to Chopin's choice that is perhaps only incompletely explained by the polar resolution across a $2p$-decacycle in the Eulerian tonnetz. We have approached these resolutions both from a triadic perspective and from a tetrachordal perspective. Why both? As we said in Section \ref{sec:Introduction}, it appears that Wagner and his contemporaries were operating within more than one harmonic system. If this seems paradoxical, keep in mind that here we work at the boundary of art and science. It is of the nature of works of art that one can never quite pin them down. They resist analysis, and yet they will readily embrace multiple modes of analysis. Each approach reveals yet another aspect of the artwork. One might even say that if the object under consideration yields completely to any one mode of analysis, then it isn't art.  In this respect, we seem to be in agreement with Tymoczko (2011) at p.\!~304: `` .\, .\, .\, I am not proposing a simple method or rule for doing musical analysis: there is no royal road to musical understanding, geometrical or otherwise.'' 

\vspace{0.20cm}

{\bf G\"otterd\"ammerung}. We turn to consider another Wagnerian example of the resolution of minor sixths to dominant sevenths, this time in an octacycle of the Tristan-genus tonnetz from {\it G\"otterd\"ammerung}. In the final section of the opera, at the beginning of the scene where Br\"unnhilde orders the vassals to pile up logs to create a great funeral pyre around the dead Siegfried's body, there is an extraordinary sequence of chords of the Tristan genus. The sequence extends over a total of thirteen measures, beginning with the $Gm^6$ at the {\it molto ritardando} just before the double bar, building up to a {\it forte} climax at the $F_\sharp M$ with the word ``Gluth,'' when the fire motive is heard, 
\begin{eqnarray}
Gm^6 \to A^7 \to B_\flat m^6 \to C^7 \to D_\flat m^6 \to E_\flat ^{7} \to Em^6 \to {C_\sharp}^{7} \to F_\sharp M.
\end{eqnarray} 
The final transition ${C_\sharp}^{7} \to F_\sharp M$  can be heard as an authentic cadence designed as a  ``wrap up,'' to break the sequence.  Looking at the Tristan-genus tonnetz, one sees that all of the tetrads in the sequence running from $Gm^6$ to $Em^6$ lie sequentially on the same octacycle. The next chord along the octacycle after the $Em^6$ is an ${F_\sharp}^{7}$. We can thus reasonably take the view that the purpose of Wagner's intervention of the ``off-cycle'' dominant seventh chord ${C_\sharp}^{7}$ just before the would-be
${F_\sharp}^{7}$ continuation of the sequence is to put the brakes on the progression with an authentic cadence, allowing it to halt at the triadic  
$F_\sharp M$ rather than continuing with a fully tetradic ${F_\sharp}^{7}$, which would demand another move to $Gm^6$ and  repetition of the cycle {\it ad infinitum} in a  {\it perpetuum mobile}.

\begin{figure}[htbp]
\includegraphics[scale=0.50]{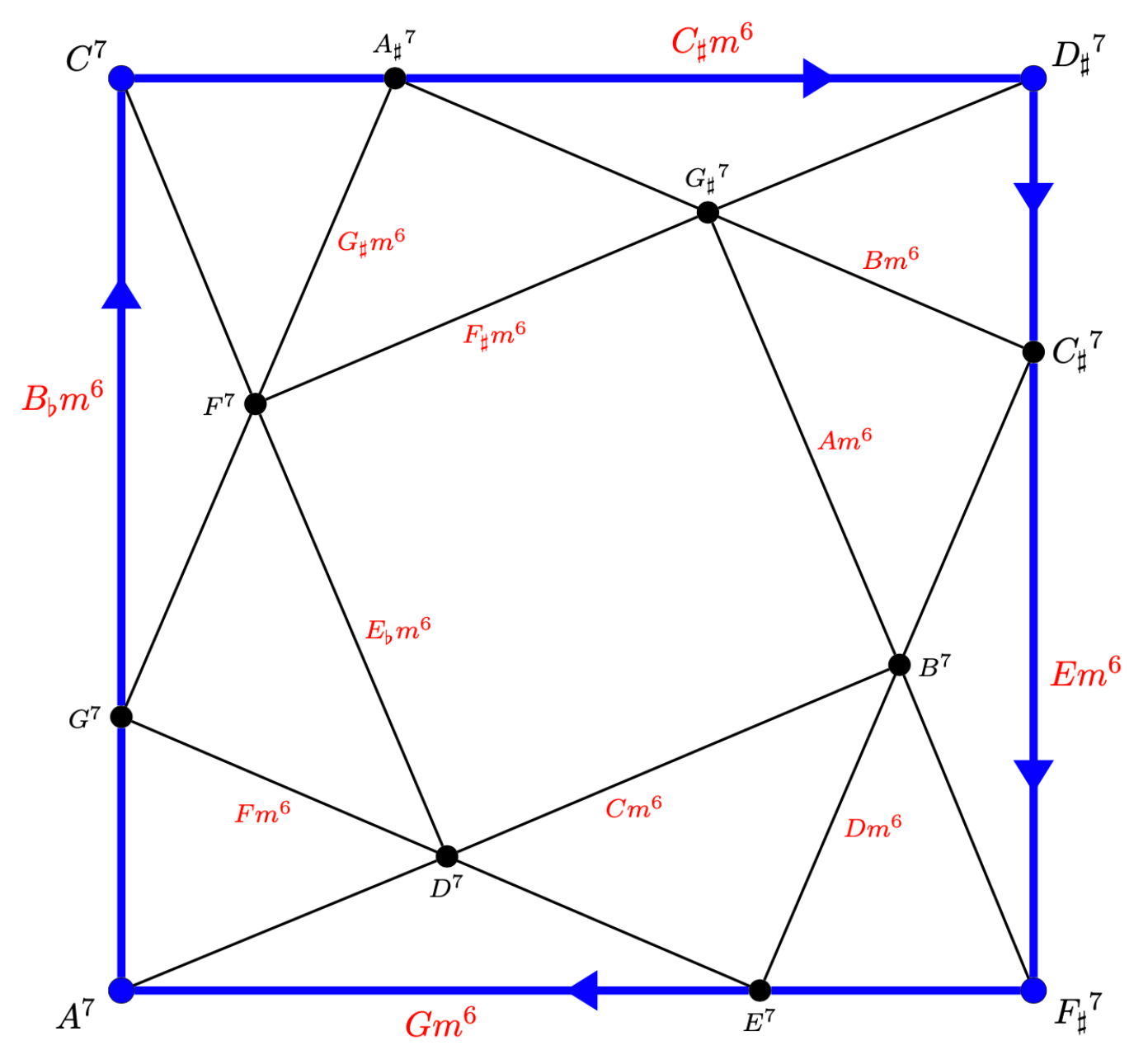}
\vspace{0cm}
\caption{In the lead up to Br\"unnhilde's final aria in {\it G\"otterd\"ammerung}, the fate motive from {\it Die Walk\"ure} is sounded in the form of a $Gm^6$ chord resolving to ${F_\sharp}^7$, announcing the intended trajectory of the progression that follows. This sequence moves from $Gm^6$ to $A^7$ to $B_\flat m^6$ to $C^7$ to $D_\flat m^6$ to ${E_\flat}^{7}$ and then to $Em^6$. Rather than resolving the $Em^6$ to a final  ${F_\sharp}^7$, Wagner breaks off the octacycle by resolving $Em^6$ to ${C_\sharp}^{7}$ which then in a standard cadence resolves to $F_\sharp M$.}
\label{fig:Brunnhilde}
\end{figure}

The geometrical structure of this progression is particularly striking if we examine the path that it makes as it transverses the configuration in Figure \ref{fig:Brunnhilde}. The context is as follows. 
At the double bar where the key changes to C major, Gutrune denounces Hagen, then turns away from the dead Siegfried and bends over the body of her brother Gunther. 
Hagen stands defiantly to one side, leaning on his spear, brooding gloomily. 

Then begins a brief orchestral interlude  in which the fate motive from {\it Die Walk\"ure} is heard in the form of a $Gm$ chord, with the tympani playing an $E$ in the bass, creating the sound of $Gm^6$,  which resolves to ${F_{\sharp} m}^7$.  This recollection of the fate motive, first heard at the Valkyrie's confrontation with Siegmund, where her destiny is determined by her fateful decision to defy Wotan, leads ultimately, after many subsequent events, to this very moment. 
This statement of the fate motive as a transition from $Gm^6$ to $F_{\sharp} m^7$ can be viewed as a  foreshadowing of the octacyclic sequence soon to follow, symbolizing the course of Br\"unnhilde's life from her days as a Valkyrie up to the present, as she prepares to enter the fire. 
The fate motive sounds again, with the fully-fledged $\varnothing^7$ harmonies of {\it G\"otterd\"ammerung}, in the form of $Am^6$ transitioning to $Bm^6$, blurred by the tympani still beating an $E$ in the bass. A third announcement of the fate motive, this time foreshortened, without the passing middle note,  begins once more at $Gm^6$, in the measure marked {\it molto ritardando}, signifying the end of the fate interlude and the onset of the octatonic sequence. 

In Figure \ref{fig:Brunnhilde}, the $Gm^6$ is represented by the horizontal line at the bottom of the configuration. 
The octatonic sequence corresponds to the quadrilateral subconfiguration running around the perimeter of the Tristan-genus configuration. The initial $Gm^6$ resolves to the $A^7$ found at the lower left-hand corner of the configuration, at the point in the score marked {\it marcato}, where there is a change of tempo, propelling the music on towards the conclusion of the music drama. After pulsing at $A^7$ there is a sense of stately movement in the next measure leading to the $B_\flat m^6$, where the music lingers, with movement again heading on to $C^7$. Half way through the pulsing of that measure, and into the next, over further movement in $C^7$, Br\"unnhilde  begins her final aria, {\it Starke Scheite schichtet mir dort}. This leads to two measures of $D_\flat m^6$, following the same pattern. 

At this point Wagner again adopts a kind of Beethovenesque foreshortening effect at the resolving ${E_\flat }^{7}$ by skipping the pulsing measure and only playing the ``movement'' measure. 
The following measure sounding an $Em^6$ is also foreshortened, with the earlier ascending dotted figures now replaced by a rapidly swirling figure, six notes to the beat, still sounding the $Em^6$ chord, suggesting something is about to change -- at which point the brakes are applied, with a swirling measure of ${C_\sharp} m^{7}$, and the terminus is reached with an $F_\sharp M$ substituting for ${F_\sharp}^{7}$, in which we hear the fire motive. Lewin (1996), pages 208-209, and Cohn (2012), pages 156-157, present analyses of this scene from a different perspective, with the sequence starting at the $A^7$. Our characterization of this important progression  of minor sixths and dominant sevenths as belonging to a single octacycle on the Tristan-genus tonnetz of Figures \ref{fig:tristan_tonnetz} and \ref{tristan_tonnetz_configuration} starting at the preceding $Gm^6$ offers what we believe may be a more satisfactory account of the matter. From our perspective, this  sequence from {\it G\"otterd\"ammerung} can be understood as a musical realization of a quadrilateral subconfiguration of the Tristan-genus tonnetz. The music thus forms another ``representation'' of the geometry.

\vspace{0.20cm}
{\bf Parsifal}. Wagner uses a version of the same progression, based on the same octacycle, in  {\it Parsifal}. Recall the sequence of four minor sixths followed by a diminished seventh that we discussed earlier, at the black knight's entrance in Act 3, shown in Figure \ref{parsifal entrance}, given by
$
B_\flat m^6 \to D_\flat m^6 \to G_\flat m^{6} \to Em^6 \to G^{\circ 7}.
$
One can easily check that the general ordering of minor sixth chords around the perimeter Hamiltonian of the Tristan-genus tonnetz is the same as the ordering of the corresponding triadic minor chords around the perimeter of the Eulerian tonnetz -- namely, ascending by fifths counterclockwise. So the conclusions that we reached earlier in our initial triadic analysis of this sequence  remain valid in the tetrachordal case, which supports our earlier assertion, in agreement with  Cohn (2012) at pages 142-145, that there is indeed some merit in a reductive approach to Wagner's tetradic harmonies. 
Specifically, we see that the initial minor sixth $B_\flat m^6$ and the final minor sixth $E m^6$ are polar on the Hamiltonian cycle in Figure  \ref{fig:tristan_tonnetz}, as they are in the triadic case, with the first resolution to the ``Tchaikovsky chord'' $D_\flat m^6$ half way between the starting point and the ending point, and the second resolution $G_\flat m^6$ at a ``staging post'' one black vertex along the cycle towards the end point. The first, second and fourth chords of the sequence lie on the octacycle containing the beginning and ending points. The final resolution to $G^{\circ 7}$, which acts again as a wrap up, is only one semitone away from the would-be continuation $G m^6$, which is another quarter turn around the octacycle. 

\begin{figure}[htbp]
\includegraphics[scale=0.60]{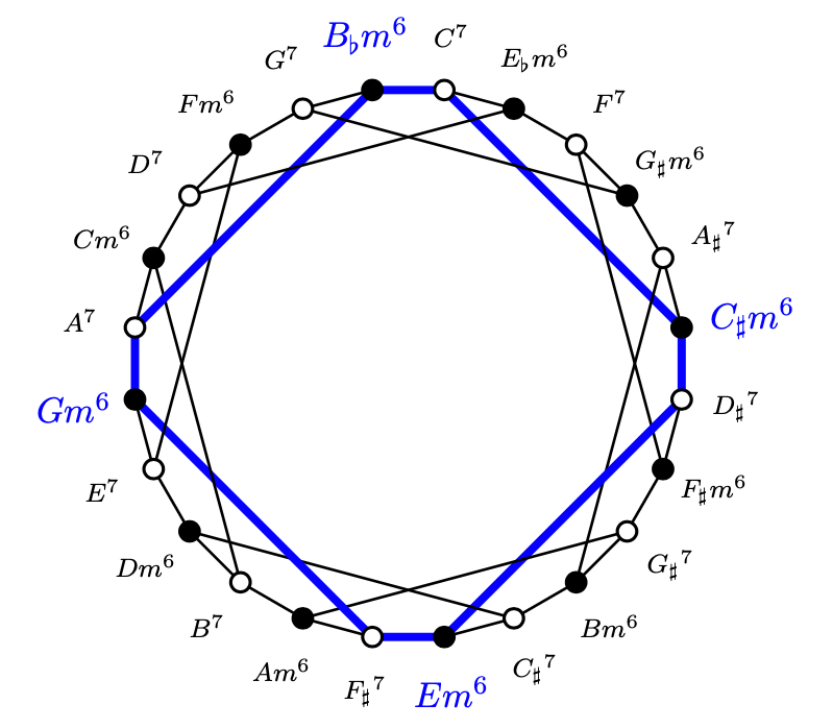}
\vspace{0cm}
\caption{Parsifal's entrance at the beginning of Act 3. The muted fanfare on $B_\flat m^6$ resolves to $D_\flat m^6$, then, after reaching a staging post at $G_\flat m^6$, moves around the octacycle to $Em^6$, before wrapping up at $G^{\circ 7}$, this diminished seventh substituting as a terminus in place of $Gm^6$. }
\label{parsifal_tetrachord_analysis}
\end{figure}

In concluding this study of configurations, tessellations and tone networks, we emphasize our view that the tonnetz and its generalizations are mathematical objects. As such, they can be represented in various ways, as networks of notes and chords, as Levi graphs, and as configurations of points and lines in Euclidean space. We have looked in detail at the Eulerian tonnetz and at a related construction that we call the Tristan-genus tonnetz which sheds light on some of the harmonies used by Wagner in his later operas. There are also links with interesting tessellations. Our general thesis remains that Levi graphs, alongside allied configurations, are fundamental as a basis for the analysis of music, for the composition of music, and even for the creation of new systems of musical composition. The Tristan-genus tonnetz provides an example of the latter, even if in that case the ``new system'' was formulated in the nineteenth century. We hope to pursue these and related topics elsewhere, including, in particular, a tonnetz for pentatonic music based on the Desargues configuration and a tonnetz for twelve-tone music based on the Cremona-Richmond configuration. 

\vspace{-0.25cm}
\begin{acknowledgments}
\noindent
LPH wishes to acknowledge E.~Chew, M.~Gotham, K.~Rietsch and other members of the Music and Acoustics Research Centre (MARC) at King's College London for stimulating discussions. The authors are grateful for comments made by participants at a MARC seminar in May 2025, and at a colloquium of the Center for Theoretical Physics of the Polish Academy of Sciences in December 2025,  where drafts of this work have been presented. We are grateful for the comments and suggestions of two anonymous reviewers. The authors also wish to thank A.~Alazemi, J.~Armstrong, D.~M.~Blasius, J.~Forth, L.~Masiero, R.~Morales, P.~Nurowski, S.~Salamon and others for helpful discussions and correspondence. 
\end{acknowledgments}

\vspace{0.30cm}

\noindent {\bf References}
%+Bibliography

\newpage

\appendix
\section*{APPENDIX: Cycle Count Tables}

\noindent We record the numbers of cycles of various types for the Eulerian tonnetz, the Archimedean tonnetz, and the Tristan-genus tonnetz. The rows are labelled by cycle length and the columns are labelled by the numbers of crossings (the $p$-numbers). Each entry in the table shows the number of cycles of a given length with the given number of crossings. Totals are tabulated for each row and each column. For example, in the case of the Eulerian tonnetz we observe that there are twelve 2$p$-hexacycles, four 3$p$-hexacycles, and a total of sixteen hexacycles altogether, and that the Eulerian tonnetz admits a grand total of 5409 cycles. 
\vspace{0.50cm}

\begin{table}[H]
\scriptsize
\centering
{\setlength{\tabcolsep}{3pt}
\begin{tabular}{|c|cccccccccccc|c|}
\hline
\multicolumn{14}{|c|}{\textbf{Cycle Count for Eulerian Tonnetz}} \\
\hline
$ $ & $0$ & $1$ & $2$ & $3$ & $4$ & $5$ & $6$ & $7$ & $8$ & $9$ & $10$ & $11$ & $\text{Total}$ \\
\hline
$6$  & $0$ & $0$ & $12$ & $4$  & $0$   & $0$   & $0$   & $0$   & $0$   & $0$   & $0$   & $0$   & $16$ \\
$8$  & $0$ & $12$& $0$  & $24$ & $3$   & $0$   & $0$   & $0$   & $0$   & $0$   & $0$   & $0$   & $39$ \\
$10$ & $0$ & $0$ & $12$ & $48$ & $60$  & $0$   & $0$   & $0$   & $0$   & $0$   & $0$   & $0$   & $120$ \\
$12$ & $0$ & $0$ & $30$ & $72$ & $78$  & $48$  & $0$   & $0$   & $0$   & $0$   & $0$   & $0$   & $228$ \\
$14$ & $0$ & $0$ & $24$ & $24$ & $228$ & $192$ & $120$ & $48$  & $0$   & $0$   & $0$   & $0$   & $636$ \\
$16$ & $0$ & $0$ & $0$  & $72$ & $168$ & $324$ & $186$ & $156$ & $0$   & $0$   & $0$   & $0$   & $906$ \\
$18$ & $0$ & $12$& $12$ & $36$ & $108$ & $324$ & $468$ & $348$ & $192$ & $0$   & $0$   & $0$   & $1500$ \\
$20$ & $0$ & $0$ & $0$  & $48$ & $42$  & $252$ & $192$ & $312$ & $174$ & $132$ & $30$  & $0$   & $1182$ \\
$22$ & $0$ & $0$ & $12$ & $0$  & $48$  & $48$  & $84$  & $108$ & $228$ & $144$ & $36$  & $12$  & $720$ \\
$24$ & $1$ & $0$ & $0$  & $0$  & $0$   & $0$   & $6$   & $36$  & $3$   & $4$   & $0$   & $12$  & $62$ \\
\hline
$\text{Total}$ & $1$ & $24$ & $102$ & $328$ & $735$ & $1188$ & $1056$ & $1008$ & $597$ & $280$ & $66$ & $24$ & $5409$ \\
\hline
\end{tabular}}
\caption{Cycle count for the Eulerian tonnetz. The rows are labelled by cycle length and the columns are labelled by the number of $p$-crossings. }
\end{table}

\begin{table}[H]
\scriptsize
\centering
{\setlength{\tabcolsep}{3pt}
\begin{tabular}{|c|cccccc|c|}
\hline
\multicolumn{8}{|c|}{\textbf{\,\,Cycle Count for Archimedean Tonnetz\,\,}} \\
\hline
$ $ & $0$ & $1$ & $2$ & $3$ & $4$ & $5$ & $\text{Total}$ \\
\hline
$4$  & $0$ & $0$ & $3$  & $0$  & $0$  & $0$  & $3$ \\
$6$  & $0$ & $6$ & $12$ & $2$  & $0$  & $0$  & $20$ \\
$8$  & $0$ & $6$ & $0$  & $18$ & $0$  & $0$  & $24$ \\
$10$ & $0$ & $0$ & $12$ & $18$ & $12$ & $6$  & $48$ \\
$12$ & $1$ & $0$ & $0$  & $2$  & $3$  & $6$  & $12$ \\
\hline
$\text{Total}$ & $1$ & $12$ & $27$ & $40$ & $15$ & $12$ & $107$ \\
\hline
\end{tabular}}
\caption{Cycle count for the Archimedean tonnetz. The rows are labelled by cycle length and the columns are labelled by the number of $p$-crossings.}
\end{table}

\begin{table}[H]
\scriptsize
\centering
{\setlength{\tabcolsep}{3pt}
\begin{tabular}{|c|cccccccccccc|c|}
\hline
\multicolumn{14}{|c|}{\textbf{Cycle Count for Tristan tonnetz}} \\
\hline
$ $ & $0$ & $1$ & $2$ & $3$ & $4$ & $5$ & $6$ & $7$ & $8$ & $9$ & $10$ & $11$ & $\text{Total}$ \\
\hline
$6$  & $0$ & $12$ & $12$ & $0$  & $0$   & $0$   & $0$   & $0$   & $0$   & $0$   & $0$   & $0$   & $24$ \\
$8$  & $0$ & $0$  & $0$  & $12$ & $3$   & $0$   & $0$   & $0$   & $0$   & $0$   & $0$   & $0$   & $15$ \\
$10$ & $0$ & $0$  & $12$ & $24$ & $60$  & $24$  & $0$   & $0$   & $0$   & $0$   & $0$   & $0$   & $120$ \\
$12$ & $0$ & $0$  & $0$  & $40$ & $78$  & $96$  & $2$   & $0$   & $0$   & $0$   & $0$   & $0$   & $216$ \\
$14$ & $0$ & $0$  & $0$  & $72$ & $120$ & $240$ & $192$ & $0$   & $0$   & $0$   & $0$   & $0$   & $624$ \\
$16$ & $0$ & $0$  & $42$ & $12$ & $144$ & $228$ & $276$ & $144$ & $0$   & $0$   & $0$   & $0$   & $846$ \\
$18$ & $0$ & $0$  & $12$ & $72$ & $108$ & $192$ & $364$ & $444$ & $192$ & $40$  & $0$   & $0$   & $1424$ \\
$20$ & $0$ & $12$ & $0$  & $24$ & $54$  & $84$  & $186$ & $372$ & $282$ & $168$ & $0$   & $0$   & $1182$ \\
$22$ & $0$ & $0$  & $12$ & $12$ & $12$  & $48$  & $96$  & $120$ & $168$ & $120$ & $24$  & $12$  & $624$ \\
$24$ & $1$ & $0$  & $0$  & $0$  & $0$   & $0$   & $0$   & $0$   & $3$   & $16$  & $18$  & $12$  & $50$ \\
\hline
$\text{Total}$ & $1$ & $24$ & $90$ & $268$ & $579$ & $912$ & $1116$ & $1080$ & $645$ & $344$ & $42$ & $24$ & $5125$ \\
\hline
\end{tabular}}
\caption{Cycle count for the Tristan tonnetz. The rows are labelled by cycle length and the columns are labelled by the number of $p$-crossings.}
\end{table}
\vspace{0.5cm}

\end{document}